\documentclass{amsart}
\usepackage{amssymb}
\usepackage[backref=page]{hyperref} 
\usepackage{pgfplots}

\numberwithin{equation}{section}

\newtheorem{theorem}{Theorem}[section]
\newtheorem{lemma}[theorem]{Lemma}
\newtheorem{proposition}[theorem]{Proposition}
\newtheorem{remark}{Remark}
\numberwithin{equation}{section}

\DeclareMathOperator{\im}{Im}
\DeclareMathOperator{\rot}{rot}
\DeclareMathOperator{\cin}{cin}
\DeclareMathOperator{\curv}{curv}

\DeclareMathOperator{\diam}{diam}
\DeclareMathOperator{\supp}{supp}

\DeclareMathOperator{\esssup}{ess\,sup}
\DeclareMathOperator{\dist}{dist}

\DeclareFontFamily{U}{mathx}{\hyphenchar\font45}
\DeclareFontShape{U}{mathx}{m}{n}{
      <5> <6> <7> <8> <9> <10>
      <10.95> <12> <14.4> <17.28> <20.74> <24.88>
      mathx10
      }{}
\DeclareSymbolFont{mathx}{U}{mathx}{m}{n}
\DeclareFontSubstitution{U}{mathx}{m}{n}
\DeclareMathAccent{\widecheck}{0}{mathx}{"71}
\DeclareMathAccent{\wideparen}{0}{mathx}{"75}

\usepackage{xcolor}
\pgfplotsset{compat=1.18}

\begin{document}

\title[Packing dimension of vertical projections in the Heisenberg group]{Packing dimension of vertical projections in the Heisenberg group}

\author{Terence L.~J.~Harris}
\address{Department of Mathematics\\ University of Wisconsin\\ 480
Lincoln Drive\\ Madison, WI 53706\\ USA}
\curraddr{School of Mathematics and Physics\\ The University of Queensland\\ St Lucia QLD 4067\\  Australia}
\email{terry.harris@uq.edu.au}

\begin{abstract}
 It is shown that if $A$ is a Borel subset of the first Heisenberg group, with Hausdorff dimension satisfying $2< \dim A < 3$, then the packing dimensions of vertical projections of $A$ are almost surely not less than $\dim A$, where both packing and Hausdorff dimensions are defined with respect to the Korányi metric. 

For the Hausdorff dimension of the projections, a weaker almost sure lower bound is obtained which improves the known bound in the range $2 < \dim A < \frac{1}{8}\left( 17 + \sqrt{33}\right) \approx 2.84$. The bound is slightly larger than $1+\frac{1}{2} \dim A$ and behaves similarly near $\dim A =2$.  Both proofs rely on a variable coefficient local smoothing inequality.
\end{abstract}

\maketitle

\section{Introduction}

Let $\mathbb{H}$ be the first Heisenberg group, identified as a set with $\mathbb{C} \times \mathbb{R} = \mathbb{R}^3$ and equipped with the product
\[ (z,t) \ast (\zeta, \tau) = \left(z+\zeta, t+\tau +\frac{1}{2} \omega(z, \zeta) \right), \]
where, for $z = x+iy$ and $\zeta = u+iv$, 
\[ \omega(z,\zeta) = -\im\left( z \overline{\zeta} \right) = z \wedge \zeta = xv - yu. \]
For each $\theta \in [0, \pi)$ let $\mathbb{V}_{\theta}^{\perp} \subseteq \mathbb{H}$ be the vertical subgroup $\{ (\lambda_1 i e^{i \theta}, \lambda_2 ) : \lambda_1, \lambda_2 \in \mathbb{R}\}$, and let $P_{\mathbb{V}_{\theta}^{\perp}} : \mathbb{H} \to \mathbb{V}_{\theta}^{\perp}$ be the vertical projection
\[ P_{\mathbb{V}_{\theta}^{\perp}}(z,t) = \left( \pi_{V_{\theta}^{\perp}}(z) , t  + \frac{1}{2}\omega( \pi_{V_{\theta}}(z), z ) \right) = (z,t) \ast P_{\mathbb{V}_{\theta}}(z,t)^{-1}, \]
where $P_{\mathbb{V}_{\theta}}: \mathbb{H} \to \mathbb{H}$ is Euclidean orthogonal projection to the line spanned by $\left( e^{i \theta}, 0\right)$, and $\pi_{V_{\theta}} : \mathbb{R}^2 \to \mathbb{R}^2$, $\pi_{V_{\theta}^{\perp}} : \mathbb{R}^2 \to \mathbb{R}^2$ are Euclidean orthogonal projection onto the span of $e^{i \theta}$, $ie^{i \theta}$ respectively.  It was conjectured in \cite[Conjecture~1.5]{BDFMT} that, if $A \subseteq \mathbb{H}$ is a Borel set, then 
\begin{equation} \label{equconjecture} \dim P_{\mathbb{V}_{\theta}^{\perp}}(A) \geq \min\{\dim A, 3 \},  \qquad \text{a.e.~$\theta \in [0, \pi)$,} \end{equation}
where $\dim$ refers to Hausdorff dimension with respect to the Korányi metric $d_{\mathbb{H}}$, given by 
\[ d_{\mathbb{H}}((z,t), (\zeta, \tau) ) = \lVert (\zeta, \tau)^{-1} \ast (z,t)  \rVert_{\mathbb{H}}, \qquad  \lVert(z,t) \rVert_{\mathbb{H}} = \left(|z|^4 + 16t^2\right)^{1/4}. \] Only the case $2 < \dim A < 3$ remains open (\cite{fasslerorponen}).  The case $\dim A \leq 1$ was solved in \cite{BDFMT}, where the problem was introduced. The previously best known bound is due to Fässler and Orponen~\cite{fasslerorponen}, who proved the conjecture \eqref{equconjecture} for $\dim A \leq 2$ and for $\dim A = 3$, and showed that for a.e.~$\theta \in [0,\pi)$,
\begin{align} \notag   \dim P_{\mathbb{V}_{\theta}^{\perp}} (A) &\geq \max\{\min\{ \dim A , 2 \}, \min\{2\dim A - 3,3 \}\} \\
\label{fo22bound} &= \begin{cases} \dim A & 0 \leq \dim A \leq 2 \\
2 & 2 < \dim A \leq 5/2 \\
2 \dim A - 3 & 5/2 < \dim A < 3 \\
3 & \dim A \geq 3. \end{cases} \end{align}
See \cite{fasslerorponen} for a brief summary of prior work on this problem. The first result of this article is the following:

\begin{theorem} \label{etccorollary1} Let $A \subseteq \mathbb{H}$ be a Borel set with $2 < \dim A < 3$. Then $\dim_P P_{\mathbb{V}_{\theta}^{\perp}}(A) \geq \dim A$ for a.e.~$\theta \in [0, \pi)$. \end{theorem}

Above, $\dim_P$ refers to the packing dimension with respect to the Korányi metric. The characterisation of packing dimension that will be used here is the upper modified box dimension:
\[ \dim_P(E) = \inf\left\{ \sup \overline{\dim}_B(E_i) : E \subseteq \bigcup_{i=1}^{\infty} E_i, \quad E_i \text{ compact} \right\}; \]
see \cite[Theorem~3.11 (g)]{cutler2} for the equivalence of this definition with the definition via packing measures, in a general metric space. The upper box dimension $\overline{\dim}_B E$ of a set $E$ is the supremum over all $s$ with the property that, for any sufficiently small $\delta>0$, there is a disjoint family of balls centred at points in $E$, of radii smaller than $\delta$, such that the sum of their radii to the power $s$ is greater than or equal to 1. By a simple pigeonholing argument, an equivalent definition results if the radii from the disjoint family of balls are required to be identical.

At some point in the proof of Theorem~\ref{etccorollary1}, the packing dimension rather than Hausdorff dimension is needed to be able to choose scales where a measure satisfies a Euclidean lower density assumption. By modifying the proof of Theorem~\ref{etccorollary1} with a recursive argument set up to find such scales, the following lower bound is obtained, which improves upon \eqref{fo22bound} in the range $2< \dim A < \frac{1}{8}\left( 17 + \sqrt{33}\right) \approx 2.84$.
\begin{theorem} \label{etccorollary2} Let $A \subseteq \mathbb{H}$ be a Borel set with $2 < t < 3$, where $t = \dim A$. Then for a.e.~$\theta \in [0, \pi)$,
\begin{align*} \dim P_{\mathbb{V}_{\theta}^{\perp}}(A) &\geq \max\left\{ 2 + \frac{(t-1)(t-2)}{t^2-4t+6}, t - \frac{ (t-2)(3-t)(t-1)}{-3t^2+14t-14},2t-3\right\} \\
\\
&= \begin{cases} 2 + \frac{(t-1)(t-2)}{t^2-4t+6} & 2 < t < \frac{5}{2} \\
t - \frac{(t-1)(t-2)(3-t)}{-3t^2+14t-14} & \frac{5}{2} \leq t < \frac{1}{8}\left( 17 + \sqrt{33}\right) \\
2t-3  &  \frac{1}{8}\left( 17 + \sqrt{33}\right) \leq t < 3. \end{cases} \end{align*} \end{theorem}

The last term $2t-3$ in the maximum is not new as it is already contained in \eqref{fo22bound}, but it does follow from the proof (though this results in an overcomplicated proof of this fact). The tangent line to $2 + \frac{(t-1)(t-2)}{t^2-4t+6}$ at $t=2$ is the function $1 + \frac{t}{2}$, which is a slightly weaker bound that follows from a less optimised version of the argument (see Remark~\ref{gammaLremark}). A comparison of \eqref{equconjecture}, Theorem~\ref{etccorollary2}, and \eqref{fo22bound} is shown in Figure~\ref{bestboundpicture}.
\begin{figure} \label{bestboundpicture}
\begin{tikzpicture}[
declare function={
	  func(\x)= (2<=\x<=3) * (\x);
		func3(\x)= (2<=\x<=3) * (max(2,2*\x-3));
		func4(\x)= (2<=\x<=2.84) * (max(2 + (\x-1)*(\x-2)/(\x*\x-4*\x+6), \x-((\x-1)*(\x-2)*(3-\x)/(-3*\x*\x+14*\x-14)));
  }]
  \begin{axis}[ 
	grid=major,
    xlabel=$\dim A$,
    ylabel={$\dim P_{\mathbb{V}_{\theta}^{\perp}}(A)$},
    ylabel style={rotate=-90},
		xtick={2,2.5,2.8430703308172535824813264335273661647775330572478490459624846838,3},
		legend style={
legend pos=outer north east,
}
  ] 
\addplot[dotted,domain=2:3]{func(x)}; 
\addplot[black,thick,domain=2:3]{func3(x)}; 
\addplot[dashed,thick,domain=2:2.8430703308172535824813264335273661647775330572478490459624846838]{func4(x)}; 
\legend{Conjecture,\cite{fasslerorponen},Theorem~\ref{etccorollary2}}
  \end{axis}
\end{tikzpicture} 
\end{figure}

A slightly simpler version of the proof of Theorem~\ref{etccorollary1} yields the following. 

\begin{proposition} \label{etccorollary} Let $\mu$ be a Borel measure on $\mathbb{H}$ which is Euclidean Ahlfors-regular. Then $\dim\left(  P_{\mathbb{V}_{\theta}^{\perp}}(\supp\mu)\right) \geq \min\left\{ 3, \dim^* \mu\right\}$ for a.e.~$\theta \in [0, \pi)$.    \end{proposition}

Even though the conclusion in Proposition~\ref{etccorollary} is for the Hausdorff dimension with respect to the Korányi metric, the Ahlfors regularity assumption is with respect to the Euclidean metric, and there is no assumption on the exponent of Euclidean Ahlfors regularity. Recall that a Borel measure $\mu$ on a metric space is Ahlfors $s$-regular if there exist positive constants $C_1$ and $C_2$ such that $C_1 r^s \leq \mu(B(x,r)) \leq C_2r^s$ for all $x \in \supp \mu$ and all $0 \leq r < \diam \supp \mu$. For any Borel measure $\mu$, the upper Hausdorff dimension of $\mu$ is defined to be
\[ \dim^*\mu = \esssup_{x \in \supp \mu} \liminf_{r \to \infty} \frac{ \log \mu(B(x,r)) }{\log r }; \]
see e.g.~(\cite[Eq.~(10.13)]{falconertechniques}) for the definition in Euclidean space. By Frostman's lemma for separable metric spaces (\cite{howroyd}), if Proposition~\ref{etccorollary} could be proved without the Euclidean Ahlfors-regular assumption, this would imply the conjectured \eqref{equconjecture}. 

\subsection{Motivation for the proof of Theorem~\ref{etccorollary1} and Proposition~\ref{etccorollary}}

The idea behind the proof of Theorem~\ref{etccorollary1} and Proposition~\ref{etccorollary} uses the Fässler-Orponen proof of the $\dim_{\mathbb{H}} A \leq 2$ case as a starting point. They prove that if $0 \leq \dim_E A \leq 1$, (where $\dim_E$ refers to Euclidean Hausdorff dimension), then $\dim_E(\pi(P_{\mathbb{V}_{\theta}^{\perp}}(A))) = \dim_E A$ for a.e.~$\theta \in [0, \pi)$, where $\pi : \mathbb{C} \times \mathbb{R} \to \mathbb{R}$ is $\pi(z,t) = t$. For $\dim_E A > 1$, it is natural to expect that $\pi(P_{\mathbb{V}_{\theta}^{\perp}}(A))$ should almost surely have positive length, but Euclidean projection theorems suggest one should expect a refinement. If $\dim_E A = s>1$, it is natural to expect that for a.e.~$\theta \in [0, \pi)$, $\pi(P_{\mathbb{V}_{\theta}^{\perp}}(A))$ should have (for any $\epsilon >0$) a positive length set of points whose fibres under $\pi \circ P_{\mathbb{V}_{\theta}^{\perp}}$ intersect $A$ in a set of Euclidean Hausdorff dimension at least $s-1-\epsilon$. A stronger refinement, which may be too strong to expect, would be that if $\dim_E A = s$, then for a.e.~$\theta \in [0, \pi)$, $\pi(P_{\mathbb{V}_{\theta}^{\perp}}(A))$ has a positive length set of points whose fibres under the restriction $\pi: \mathbb{V}_{\theta}^{\perp} \to \mathbb{R}$ intersect $P_{\mathbb{V}_{\theta}^{\perp}}(A)$ in a set of Euclidean Hausdorff dimension at least $s-1-\epsilon$. If this stronger refinement were true, then a simple Fubini-type argument (see \eqref{claimedinequality} below) with Euclidean-Korányi dimension comparison would yield the conjectured inequality~\eqref{equconjecture} for Korányi-Hausdorff dimension.  However, a discrete counterexample of Orponen from 2022 \cite{orponenprivate} suggests that $s-1$ is not possible above when $1 \leq \dim_E A \leq 2$, and the best one could hope for is probably $(s-1)/2$, at least for a discretised analogue of the problem. For this reason, the Korányi Hausdorff dimension $\dim A$ is used below in the domain to avoid the Euclidean-Korányi dimension comparison step. 

For $2 < s \leq 3$, let $\beta(s)$ be supremum over all $\beta \geq 0$ with the property that, for any Borel set $A \subset \mathbb{H}$ with $\dim_{\mathbb{H}} A = s$, for a.e.~$\theta \in [0, \pi)$, the set $\pi(P_{\mathbb{V}_{\theta}^{\perp}}(A))$ has a positive length set of points whose fibres under the restriction $\pi: \mathbb{V}_{\theta}^{\perp} \to \mathbb{R}$ intersect $P_{\mathbb{V}_{\theta}^{\perp}}(A)$ in Hausdorff dimension at least $\beta$. It seems reasonable to conjecture that $\beta(s) \geq s-2$ for $2 < s \leq 3$. If this were true, then a simple Fubini-type dimension comparison argument would yield \eqref{equconjecture}. Theorem~\ref{intersectiontheorem2} shows that this conjecture is true if the Hausdorff dimension of the fibres is replaced by packing dimension. Again, a simple Fubini-type argument yields Theorem~\ref{etccorollary1} as a corollary.

A version of this problem for measures is, given a Borel measure $\mu$ supported in the unit ball of $\mathbb{H}$ satisfying an $s$-dimensional Frostman condition with respect to the Korányi metric, is it true that for any $\epsilon >0$, for a.e.~$\theta \in [0, \pi)$, there is a positive length set of points in the vertical axis whose fibres under the restriction $\pi: \mathbb{V}_{\theta}^{\perp} \to \mathbb{R}$ intersect the support of $P_{\mathbb{V}_{\theta \sharp}^{\perp}}\mu$ in dimension at least $s-2-\epsilon$? This is proved in Theorem~\ref{intersectiontheorem2} under the assumption that $\mu$ is Euclidean Ahlfors-regular. Similarly to the above, a simple Fubini-type argument yields Proposition~\ref{etccorollary} as a corollary. Proposition~\ref{etccorollary} was proved before Theorem~\ref{etccorollary1}, but then it was noticed that the non-dependence on the exponent of Euclidean Ahlfors-regularity meant that the idea extends to prove Theorem~\ref{etccorollary1}. The reason that packing rather than Hausdorff dimension works is that the packing dimension allows the selection of a potentially sparse sequence of scales at which $\mu$ looks like a Euclidean ``semi Ahlfors-regular'' measure at these scales.

The (probably sharp) projection theorem for $P_{\mathbb{V}_{\theta}^{\perp}}$ with Euclidean metric in domain and co-domain is $\dim_E P_{\mathbb{V}_{\theta}^{\perp}}(A) \geq (1+\dim_E A)/2$ when $1 \leq \dim_E A \leq 2$. This was originally proved by S.~Wu in 2024, but not published. The (conjectured) sharpness of this bound is related to the discrete counterexample of Orponen from 2022 \cite{orponenprivate} mentioned above. 

An important tool in the proof is a Euclidean $L^{p}$ inequality for projections $\pi \circ P_{\mathbb{V}_{\theta}^{\perp}}$ proved in Section~\ref{lpestimates}. The setup of the argument to convert this into an intersection theorem borrows from the method in \cite{mattila3}, to convert $L^p$ inequalities for projections into results about intersections. 

To prove the $L^p$ inequality for projections $\pi \circ P_{\mathbb{V}_{\theta}^{\perp}}$ in Section~\ref{lpestimates}, a duality idea, based on the point-curve duality from \cite{fasslerorponen}, is used in Lemma~\ref{adjoint} to convert it into an inequality for an averaging operator over curves, which is deduced from the variable coefficient local smoothing inequality of Gao-Liu-Miao-Xi~\cite{gaoliumiaoxi}. The local smoothing inequality of Beltran-Hickman-Sogge~\cite{beltranhickmansogge}, which holds for a more restricted range of exponents, would be just as useful for the application here, as the inequality is only needed for some finite exponent. The local smoothing inequality from~\cite{gaoliumiaoxi} is a variable coefficient version of the local smoothing inequality for the wave equation in $\mathbb{R}^{2+1}$ of Guth-Wang-Zhang~\cite{guthwangzhang}. Some of the Kakeya-type inequalities from \cite{guthwangzhang} were used in \cite{fasslerorponen} to prove the $\dim_{\mathbb{H}} A = 3$ case of the vertical projection problem, but the application of local smoothing here is very different to that in \cite{fasslerorponen}.

 The proof of the $L^p$ inequality for projections $\pi \circ P_{\mathbb{V}_{\theta}^{\perp}}$ in Section~\ref{lpestimates} is inspired by the proof of \cite[Corollary 3]{wolff}, but a direct imitation of the proof of Corollary~3 in \cite{wolff} would only yield positive length of projections $\pi \circ P_{\mathbb{V}_{\theta}^{\perp}}$, and a bit more care is needed to obtain an $L^p$ bound with $p>1$. 

An important ingredient for proving the intersection theorem is a quantitative projection lemma for vertical projections with Korányi metric in the domain and Euclidean metric in co-domain, given in Lemma~\ref{projectionconjecture} below. This lemma is not precise enough to deduce Theorem~\ref{etccorollary1} or Theorem~\ref{etccorollary2} directly, as it gives no information about the ``quasi-product'' structure of the projections, but it will be important in a sub-case of the proof of Theorem~\ref{etccorollary1} and Theorem~\ref{etccorollary2}. In Section~\ref{lightraykakeyasection}, Lemma~\ref{projectionconjecture} is deduced as a corollary of the $L^{3/2}$ bound on projections from \cite{harris2023}, which in turn used many of the ideas from \cite{fasslerorponen}. The use of the $L^{3/2}$ bound from \cite{harris2023} could probably be substituted by the $L^2$ bound from \cite{fasslerorponen}. Moreover, the use of the $L^{3/2}$ bound from \cite{harris2023} could also be substituted by a slightly weaker $L^{3/2}$ bound allowing $\delta^{-\epsilon}$ losses, which would permit a simpler proof using the non-endpoint trilinear Kakeya inequality in place of the endpoint version (see \cite{guthshortproof}). 

\subsection{Rough description of the proof of the Hausdorff dimension bound (Theorem~\ref{etccorollary2})}
In the proof of the Hausdorff dimension lower bound, initially following the proof of the Euclidean Ahlfors-regular case, one needs to estimate some integral $I(\Delta_L)$ depending on a parameter $\Delta_L$. If a lower density assumption is satisfied as in the Euclidean-Ahlfors regular case, then a good bound on $I(\Delta_L)$ is obtained. Otherwise, one can trivially bound $I(\Delta_L) \leq I(\Delta_{L-1})$ for some courser scale $\Delta_{L-1} > \Delta_L$. This step looks very blunt and may cause loss, but the failure of the lower density assumption at scale $\Delta_L$ implies a gain in some higher frequency terms in the bound of $I(\Delta_{L-1})$. If a lower density assumption at scale $\Delta_{L-1}$ is satisfied then a good bound is obtained on $I(\Delta_{L-1})$, and otherwise $\Delta_{L-1}$ can be replaced by a courser scale $\Delta_{L-2}$ and this process continues. When $\dim A>2$ this recursive process terminates soon enough that a gain over the bound of 2 is possible.

\section{A quantitative projection lemma with Korányi metric in domain and Euclidean metric in co-domain} \label{lightraykakeyasection}

Given a measure $\mu$ on a measurable space $(X, \mathcal{A})$, and measurable function $f: X \to Y$ from $X$ into a measurable space $(Y, \mathcal{B})$, the pushforward $f_{\sharp}\mu$ of $\mu$ under $f$ is defined by $\left(f_{\sharp}\mu\right)(E) = \mu(f^{-1}(E))$ for any $E \in \mathcal{B}$. Equivalently, for any non-negative measurable function $g$ on $Y$, $\int g \, d\left( f_{\sharp}\mu \right) = \int (g \circ f) \, d\mu$.  The pushforward is defined similarly for complex measures. 
For a Borel measure $\nu$ on $\mathbb{H}$, define 
\[ c_{t,\mathbb{H}}(\nu)  = \sup_{x \in \mathbb{H}, r >0 } r^{-t} \nu\left( B_{\mathbb{H}}(x,r) \right). \]

The following lemma is not used for the packing dimension result Theorem~\ref{etccorollary1}, but it is needed for the Hausdorff dimension bound Theorem~\ref{etccorollary2}. Roughly, it states that Euclidean mollification of a measure does not increase its dimension with respect to the Korányi metric. 

\begin{lemma} \label{euclidean} Let $t \in [0, 4]$. 
Let $\phi$ be a smooth bump function equal to 1 on $B_E(0,1)$ and vanishing outside $B_E(0,2)$, and let $\phi_{\delta}(x) = \delta^{-3} \phi(x/\delta)$. Let $\mu$ be a Borel measure on $B_{\mathbb{H}}(0,10)$. 
Then 
\[ c_{t, \mathbb{H}}( \mu \ast_E \phi_{\delta}) \lesssim c_{t, \mathbb{H}}(\mu), \]
where $\ast_E$ denotes Euclidean convolution and $\ast_{\mathbb{H}}$ denotes convolution in $\mathbb{H}$. \end{lemma}
\begin{proof}  Let $\mathcal{B}$ be a boundedly overlapping cover of $\mathbb{H}$ by Korányi $\delta$-balls. Then
\[ \mu \ast_E \phi_{\delta} \lesssim \nu \ast_E \phi_{100\delta}, \]
where  
\[ \nu =  \sum_{B \in \mathcal{B}} \mu(B) \delta^{-4} \chi_{2B}, \] 
so it suffices to prove that $c_{t, \mathbb{H}}(\nu \ast_E \phi_{100\delta}) \lesssim c_{t, \mathbb{H}}(\mu)$. For each $k$, let $\mathcal{B}_k$ consist of the Korányi $\delta$-balls from $\mathcal{B}$ with $2^{-(k+1)} \leq \mu(B) < 2^{-k}$. 
For each $k$, let 
\[ \nu_k  =\sum_{B \in \mathcal{B}_k} \mu(B) \delta^{-4} \chi_{B}, \qquad \nu_{k,C} = \sum_{B \in \mathcal{B}_k} \mu(B) \delta^{-4} \chi_{CB}, \]
for some large constant $C$ to be chosen. 
Let $Z(\mathbb{H}) = \{(0,0)\} \times \mathbb{R}$ be the centre of $\mathbb{H}$. Let $\{x_1, \dotsc, x_M\}$ be a maximal subset of $B_{E}(0, 10^4 \delta) \cap Z(\mathbb{H})$ which is $\delta$-separated in the Korányi metric. Then: 
\begin{equation} \label{coveringproperty} \{ B_{\mathbb{H}}(x_j, 10^{5} \delta)\}_{j=1}^M \text{ covers } B_E(0,10^4\delta), \end{equation}
 and $M \sim \delta^{-1}$. It will be shown that for each $k$,
\begin{equation} \label{enough} \nu_k \ast_E \phi_{100\delta} \lesssim \delta \sum_{j=1}^M L_{x_j \sharp} \nu_{k,C}, \end{equation}
where $L_{x_j}(y) = x_j \ast y$ is left translation by $x_j$ (where $\ast$ is the Heisenberg group product). It is first shown that \eqref{enough} implies the lemma. If \eqref{enough} holds, then 
\[ \nu \ast_E \phi_{100\delta} \leq \sum_k \nu_k \ast_E \phi_{100\delta} \lesssim  \delta \sum_k  \sum_{j=1}^M L_{x_j \sharp} \nu_{k,C} = \delta \sum_{j=1}^M L_{x_j \sharp}\left( \sum_k    \nu_{k,C}\right). \]
By denoting $\nu_C = \sum_{B \in \mathcal{B}} \mu(B) \delta^{-4} \chi_{CB}$, this gives 
\[ \nu \ast \phi_{100\delta} \lesssim  \delta \sum_{j=1}^M L_{x_j \sharp}(\nu_C). \]
Hence, by left invariance and since $M \sim \delta^{-1}$, 
\[ c_{t, \mathbb{H}}\left( \nu \ast \phi_{\delta}\right) \lesssim \delta \sum_{j=1}^M c_{t, \mathbb{H}}\left( L_{x_j \sharp} \nu_C \right)  = \delta \sum_{j=1}^M c_{t, \mathbb{H}}( \nu_C) \lesssim c_{t, \mathbb{H}}\left( \mu\right), \]  
This shows that \eqref{enough} implies the lemma, so it remains to establish \eqref{enough}. Let $x \in \mathbb{H}$, and let 
\[ N = \left\lvert\left\{ B \in \mathcal{B}_k : \dist_E(x, B) < 200 \delta \right\}\right\rvert. \] 
Then 
\begin{equation} \label{firsthand} (\nu_k \ast \phi_{100\delta})(x) \lesssim 2^{-k} \delta^{-3} N . \end{equation}
On the other hand, 
\begin{multline*} \delta \sum_{j=1}^M L_{x_j \sharp} \nu_{k,C}(x)  = \delta \sum_{j=1}^M \sum_{B \in \mathcal{B}_k} \mu(B) \delta^{-4} \chi_{L_{x_j}(CB)}(x) \\
\sim 2^{-k}\delta^{-3} \sum_{B \in \mathcal{B}_k} \left\lvert \left\{ j \in \{1, \dotsc, M\} : x \in L_{x_j}(CB) \right\} \right\rvert \\
\geq 2^{-k}\delta^{-3} \sum_{B \in \mathcal{B}_k : \dist_E(x,B)<100\delta} \left\lvert \left\{ j \in \{1, \dotsc, M\} : x \in L_{x_j}(CB) \right\} \right\rvert. \end{multline*} 
Since the sum is over $N$ terms, by comparing with \eqref{firsthand} it suffices to show that for each $B \in \mathcal{B}$ with $\dist_E(x,B) < 100\delta$, there is at least one $j \in \{1, \dotsc, M\}$ with $x \in L_{x_j}(CB)$. By denoting the centre of $B$ by $x_B$, it suffices to prove that there is at least one $j \in \{1, \dotsc, M\}$ with $d_{\mathbb{H}}(x, x_j \ast x_B) < C\delta$. But since $x_j \in Z(\mathbb{H})$ for all $j$, this is equivalent to finding at least one $j$ such that $d_{\mathbb{H}}(x,  x_B \ast x_j) < C\delta$, or equivalently $d_{\mathbb{H}}(x_B^{-1}\ast x,  x_j) < C\delta$ by left-invariance. But $x_B^{-1} \ast x \in B_E(0, 10^4 \delta)$ by local Lipschitz continuity of $F(y) = x_B^{-1} \ast y$ with respect to the Euclidean metric and since $d_E(x, x_B) < 200\delta$, and therefore there must exist $j \in \{1, \dotsc, M\}$ with $d_{\mathbb{H}}(x_B^{-1} \ast x, x_j) < 10^{5}\delta$ by the covering property \eqref{coveringproperty}. Taking $C > 10^{5}$ finishes the proof.   \end{proof} 

The following theorem converts the $L^{3/2}$ projection bound from~\cite{harris2023} into a quantitative projection theorem for the vertical projections, with Euclidean metric in the co-domain and Korányi metric in the domain.

\begin{theorem} \label{projectionconjecture} Suppose that $2 \leq t \leq 3$, and that $\nu$ is a Borel measure supported in the unit ball of $\mathbb{H}$ such that $c_{t,\mathbb{H}}(\nu)  < \infty$. Then, for any $\epsilon >0$, there exists $\delta_0 >0$ and a sufficiently small $\eta>0$ depending only on $t$ and $\epsilon$, such that for all $0 < \delta \leq \delta_0$, 
\begin{multline} \label{zset} \nu\left\{ x \in \mathbb{H} : \mathcal{H}^1\left\{ \theta \in [0, \pi) : P_{\mathbb{V}_{\theta \sharp}^{\perp}} \nu\left( B_E \left( P_{\mathbb{V}_{\theta}^{\perp}}(x), \delta\right) \right)  \geq  c_{t,\mathbb{H}}(\nu) \delta^{t-1 - \epsilon } \right\} \geq \delta^{\eta} \right\} \\
\leq \nu(\mathbb{H}) \delta^{\eta}. \end{multline} 
More generally, if 
\[ C(t, \nu) := \sup_{r>0} \min\left\{ r \left\lVert \nu \ast_E \phi_{\delta} \right\rVert_{\infty}, r^{t-3} c_{t, \mathbb{H}}(\nu) \right\}, \]
where $\phi_{\delta}(x)= \delta^{-3} \phi(x/\delta)$ for $\phi$ a smooth bump function equal to 1 on $B_E(0,1)$ and vanishing outside $B_E(0,2)$, then \eqref{zset} can be replaced by  
\begin{multline} \label{zset2} \nu\left\{ x \in \mathbb{H} : \mathcal{H}^1\left\{ \theta \in [0, \pi) : P_{\mathbb{V}_{\theta \sharp}^{\perp}} \nu\left( B_E \left( P_{\mathbb{V}_{\theta}^{\perp}}(x), \delta\right) \right)  \geq  C(t, \nu) \delta^{2 - \epsilon } \right\} \geq \delta^{\eta} \right\} \\
\leq \nu(\mathbb{H}) \delta^{\eta}. \end{multline} 
 \end{theorem}

\begin{remark} Lemma~\ref{projectionconjecture} can roughly be interpreted as saying that, for a typical point $x$ in the support of $\nu$, the pushforward measure of $\nu$ under vertical projection for a typical $\theta$ satisfies a $(t-1)$-dimensional Frostman condition on the Euclidean $\delta$-disc whose inverse under $P_{\mathbb{V}_{\theta}^{\perp}}$ is the (horizontal or $SL_2$) $\delta$-tube through $x$. This kind of formulation of a projection theorem (for a different family of projections) first appeared in \cite{orponenvenieri}.  \end{remark}

\begin{remark} The simpler version \eqref{zset} is enough for the packing dimension bound (Theorem~\ref{etccorollary1}), but the more detailed version \eqref{zset2} is needed for the Hausdorff dimension bound (Theorem~\ref{etccorollary2}). Using only \eqref{zset} in the Hausdorff dimension argument would give only the weaker lower bound $2+\frac{1}{2} \left( \dim A - 2\right)^2$ instead of $1 + \frac{1}{2} \dim A$.  \end{remark}

\begin{proof}[Proof of Lemma~\ref{projectionconjecture}] By the dimension comparison principle (see \eqref{comparisonfrostman} below), 
\[ \|\nu \ast_E \phi_{\delta} \|_{\infty} \lesssim c_{t, \mathbb{H}}(\nu) \delta^{t-4}. \]
Hence 
\[ C(t, \nu) \lesssim c_{t, \mathbb{H}}(\nu)\sup_{r>0}\min\left\{ r\delta^{t-4}, r^{t-3}\right\} = c_{t,\mathbb{H}}(\nu) \delta^{t-3}. \]
Substituting into \eqref{zset2} shows that \eqref{zset2} implies \eqref{zset}, so it suffices to prove \eqref{zset2}.

Let $\mu = \nu \ast_{E} \eta_{\delta}$, where $\eta_{\delta}(x) = \delta^{-3}\eta(x/\delta)$, with $\eta$ a non-negative smooth bump function supported in $B_{E}(0,1)$, such that $\eta \sim 1$ on $B_{E}(0,1/2)$ and $\int_{\mathbb{H}} \eta \, d\mathcal{H}^3_E =1$.  Since the projections $P_{\mathbb{V}_{\theta}^{\perp}}$ are locally Lipschitz when considered as functions from $\left( \mathbb{H}, d_{E}\right)$ to $\left( \mathbb{V}_{\theta}^{\perp}, d_E \right)$, for any $x \in B_{\mathbb{H}}(0,1)$ and $y \in \mathbb{H}$ with $d_{E}(x,y) < \delta$, and any $\theta \in [0, \pi)$,
\[ P_{\mathbb{V}_{\theta \sharp}^{\perp}} \mu\left( B_E \left( P_{\mathbb{V}_{\theta}^{\perp}}(y), 100\delta\right) \right)  \gtrsim P_{\mathbb{V}_{\theta \sharp}^{\perp}} \nu\left( B_E \left( P_{\mathbb{V}_{\theta}^{\perp}}(x), \delta\right) \right); \]
by unpacking the definitions in the left-hand side and applying Fubini. Therefore, if $Z'$ is the set from \eqref{zset}:
\[ Z' = \left\{ x \in \mathbb{H} : \mathcal{H}^1\left\{ \theta \in [0, \pi) : P_{\mathbb{V}_{\theta \sharp}^{\perp}} \nu\left( B_E \left( P_{\mathbb{V}_{\theta}^{\perp}}(x), \delta\right) \right)  \geq   C(t, \nu) \delta^{2 - \epsilon } \right\} \geq \delta^{\eta} \right\}, \]
then taking a maximal $\sim \delta$-separated subset of $Z'$ in the Euclidean metric to get a boundedly overlapping cover of $Z'$ by $\sim \delta$ Euclidean balls $B$, using that $\nu(B) \lesssim \mu(B)$, letting 
\begin{multline*} Z = \\
\left\{ x \in \mathbb{H} : \mathcal{H}^1\left\{ \theta \in [0, \pi) : P_{\mathbb{V}_{\theta \sharp}^{\perp}} \mu\left( B_E \left( P_{\mathbb{V}_{\theta}^{\perp}}(x), 100\delta\right) \right)  \gtrsim  C(t, \nu) \delta^{2 - \epsilon } \right\} \geq \delta^{\eta} \right\}, \end{multline*}
and using that $\bigcup B \subseteq Z$,  yields
\[ \nu(Z') \lesssim \mu(Z). \]
Therefore, using $\mu(\mathbb{H}) \lesssim \nu(\mathbb{H})$, it suffices to show that $\mu(Z) \leq \delta^{2\eta}\mu(\mathbb{H})$. 

 Let $p= 3/2$. By two applications of Chebychev's inequality, 
\begin{multline*} \mu(Z) \lesssim  \\\delta^{-(2-\epsilon)(p-1)-\eta} C(t, \nu)^{-(p-1)} \int \int_0^{\pi} \left(P_{\mathbb{V}_{\theta \sharp}^{\perp}} \mu\left( B_E \left( P_{\mathbb{V}_{\theta}^{\perp}}(x), 100\delta\right) \right) \right)^{p-1} \, d\theta \, d\mu(x). \end{multline*}
Using Fubini and the definition of pushforward, this can be simplified to
\begin{multline*} \mu(Z) \lesssim \delta^{-(2-\epsilon)(p-1)-\eta} C(t, \nu)^{-(p-1)} \times  \\ \int_0^{\pi} \int  \left(P_{\mathbb{V}_{\theta \sharp}^{\perp}} \mu\left( B_E \left( x, 100\delta\right) \right) \right)^{p-1} \, d\left( P_{\mathbb{V}_{\theta \sharp}^{\perp}} \mu\right)(x) \, d\theta . \end{multline*}
This can be written as
\begin{multline*} \mu(Z) \lesssim \delta^{\epsilon(p-1)-\eta} C(t, \nu)^{-(p-1)} \times  \\\int_0^{\pi} \int  \left(\delta^{-2} P_{\mathbb{V}_{\theta \sharp}^{\perp}} \mu\left( B_E \left( x, 100\delta\right) \right) \right)^{p-1} \, d\left( P_{\mathbb{V}_{\theta \sharp}^{\perp}} \mu\right)(x) \, d\theta. \end{multline*} 
If $M_{\theta}$ is the Hardy-Littlewood maximal operator on $L^{3/2}(\mathbb{V}_{\theta}^{\perp})$ (identified with $L^{3/2}(\mathbb{R}^2)$), the above gives
\[\mu(Z) \lesssim  \delta^{\epsilon(p-1)-\eta} C(t, \nu)^{-(p-1)} \int_0^{\pi} \int_{\mathbb{V}_{\theta}^{\perp}} \left\lvert M_{\theta} P_{\mathbb{V}_{\theta \sharp}^{\perp}} \mu\right\rvert^{p} d\mathcal{H}^2_E \, d\theta,\]
where $\mathcal{H}^2_E$ is the area or Lebesgue measure on $\mathbb{V}_{\theta}^{\perp}$. By the boundedness of the Hardy-Littlewood maximal operator on $L^{3/2}(\mathbb{R}^2)$, applied to each $\theta$, the above gives
\[ \mu(Z) \lesssim  \delta^{\epsilon(p-1)-\eta} C(t, \nu)^{-(p-1)} \int_0^{\pi} \int_{\mathbb{V}_{\theta}^{\perp}} \left\lvert P_{\mathbb{V}_{\theta \sharp}^{\perp}} \mu\right\rvert^{p} d\mathcal{H}^2_E \, d\theta. \]  
By \cite[Theorem~3.1]{harris2023}, which has $p=3/2$, this gives 
\begin{equation} \label{pause1000} \mu(Z) \lesssim  \delta^{\epsilon(p-1)-\eta} C(t, \nu)^{-(p-1)} c_{3+\epsilon^2,\mathbb{H}}(\mu)^{p-1} \mu(\mathbb{H}), \end{equation}
where the implicit constant is allowed to depend on $\epsilon$. 
For any $r>0$, by Lemma~\ref{euclidean}, 
\[ r^{-(3+\epsilon^2)}\mu(B_{\mathbb{H}}(x,r)) \lesssim r^{-\epsilon^2} \min\left\{ r \left\lVert \nu \ast_E \phi_{\delta} \right\rVert_{\infty}, r^{t-3} c_{t, \mathbb{H}}(\nu) \right\}. \]
This yields
\begin{equation} \label{updown} c_{3+\epsilon^2, \mathbb{H}}(\mu) \lesssim \sup_{r>0}  r^{-\epsilon^2} \min\left\{ r \left\lVert \nu \ast_E \phi_{\delta} \right\rVert_{\infty}, r^{t-3} c_{t, \mathbb{H}}(\nu) \right\} \lesssim \delta^{-\epsilon^2} C(\nu, t). \end{equation}
In the minimum, the first function is increasing in $r$ and the second function is decreasing in $r$, and the above used that $\left\lVert \nu \ast_E \phi_{\delta} \right\rVert_{\infty} \lesssim \delta^{-3}r^{t-1} c_{t, \mathbb{H}}(\nu)$ (from \eqref{comparisonfrostman} below) to verify that the sup above occurs when $r \gtrsim \delta$. Substituting \eqref{updown} into \eqref{pause1000} gives $\mu(Z) \lesssim \delta^{(\epsilon-\epsilon^2)(p-1)-\eta} \mu(\mathbb{H})$. Taking $\eta = \epsilon/100$ gives $\mu(Z) \leq \delta^{2\eta} \mu(\mathbb{H})$ for $\delta$ sufficiently small, and by the reasoning above, this finishes the proof.   \end{proof}

\section{An \texorpdfstring{$L^p$}{Lp} inequality for vertical projections in the Euclidean metric} \label{lpestimates}
Recall that $\pi: \mathbb{H} \to \mathbb{R}$ is the projection $(z,t) \mapsto t$ onto the vertical axis (identified with $\mathbb{R}$).

\begin{lemma} \label{adjoint} The formal adjoint of the ``rotating projection'' operator $T$ defined by 
\[ Tf(\theta, r) =  \left(\pi_{\sharp } P_{\mathbb{V}_{\theta \sharp}^{\perp} }f\right)(r) \]
 is the averaging operator $A$ defined by 
\[ Ag(z,t) = \int_0^{\pi} g\left(\theta, t + \frac{1}{2} \omega(\pi_{V_{\theta}}(z), z ) \right) \, d\theta, \]
where $z \in \mathbb{R}^2$ and $t \in \mathbb{R}$. More precisely, if $f$ is in $C_0^{\infty}(\mathbb{R}^3)$ (identified with a measure) and $g \in C_0^{\infty}([0, \pi] \times \mathbb{R})$, then 
\[  \int_0^{\pi} \int_{\mathbb{R}} Tf(\theta,r)  g(\theta, r) \, dr \, d\theta = \int_{\mathbb{R}^3} f(z,t) Ag(z,t) \, dz \, dt. \]  \end{lemma}
\begin{proof} For each $\theta \in [0, \pi]$, by the definition or characterisation of pushforward measures, 
 \begin{multline*} \int_{\mathbb{R}} Tf(\theta,r)  g(\theta, r) \, dr = \int  g(\theta, r) \, d\left(\pi_{\sharp } P_{\mathbb{V}_{\theta \sharp}^{\perp} }f\right)(r) \\
= \int_{\mathbb{R}^3} f(z,t) g\left( \theta, \pi\left( P_{\mathbb{V}_{\theta}^{\perp}}(z,t)\right) \right) \, dz \, dt. \end{multline*} 
Integrating in $\theta$, using the formula $P_{\mathbb{V}_{\theta}^{\perp}}(z,t) = \left(\pi_{V_{\theta}}(z), t + \frac{1}{2} \omega\left( \pi_{V_{\theta}}(z), z \right) \right)$, and then Fubini, gives
\[ \int_0^{\pi} \int_{\mathbb{R}} Tf(\theta,r)  g(\theta, r) \, dr \, d\theta = \int_{\mathbb{R}^3} \left[\int_0^{\pi}g\left(\theta, t + \frac{1}{2} \omega( \pi_{V_{\theta}}(z), z) \right) \, d\theta\right] f(z,t) \, dz \, dt.  \]
This proves the lemma.
\end{proof}

In the theorem below, $c_{\alpha}(\mu) = c_{\alpha, E}(\mu)$ is defined with respect to the Euclidean metric, i.e.~ $c_{\alpha,E}(\mu) = \sup_{x \in \mathbb{H}, r >0} \frac{ \mu(B_E(x,r) ) }{r^{\alpha} }$. 

\begin{theorem} \label{lpestimate} Let $\alpha >1$ and $1 < p \leq 4/3$. Then for any $\epsilon >0$, the following holds for all $R \geq 1$. 

If $\mu$ is a Borel measure supported in a Euclidean ball of radius $R^{-1}$, such that $|z| \sim 1$ for all $(z,t)$ in the support of $\mu$, with $c_{\alpha,E}(\mu)< \infty$, then
\begin{equation} \label{firstestimate} \int_0^{\pi} \int_{\mathbb{R} } \left\lvert \pi_{\sharp } P_{\mathbb{V}_{\theta \sharp}^{\perp} }\mu \right\rvert^p \, d\mathcal{H}^1_E \, d\theta \leq C_{\alpha, \epsilon}  c_{\alpha,E}(\mu)^{p-1}  \mu(\mathbb{H}) R^{\epsilon -(\alpha-1)(p-1)}, \end{equation}
 In particular, $\pi_{\sharp } P_{\mathbb{V}_{\theta \sharp}^{\perp} }\mu \ll \mathcal{H}^1_E$ for a.e.~$\theta \in [0, \pi)$ whenever $\alpha >1$ and $\mu$ is a compactly supported Borel measure satisfying the Euclidean Frostman condition $c_{\alpha,E}(\mu)< \infty$.  

Suppose that the assumption that $\mu$ is supported in a Euclidean ball of radius $R^{-1}$ is replaced by the assumption that $\mu$ is supported in a Euclidean ball of radius $\sim 1$, still with $|z| \sim 1$ for all $(z,t)$ in the support of $\mu$. If $\psi_R$ is a smooth bump function on $\lvert \xi \rvert \gtrsim R$,  then
\begin{multline} \label{thirdestimate} \left(\int_0^{\pi} \int_{\mathbb{R} } \left\lvert \pi_{\sharp } P_{\mathbb{V}_{\theta \sharp}^{\perp} }\left(\mu \ast \widecheck{\psi_R} \right) \right\rvert^p \, d\mathcal{H}^1_E \, d\theta\right)^{1/p} \\
\leq C_{ \epsilon} \mu(\mathbb{H})^{1/p} \sum_{j \geq \log_2 R}   2^{j\epsilon} \left(\sup_{x \in \mathbb{H}} 2^j \mu(B_E(x,2^{-j}))\right)^{1/p'},  \end{multline}
where the convolution is Euclidean. In particular, if $C$ is a constant such that $\mu(B_E(x,r)) \leq Cr^{\alpha}$ for all $x \in \supp \mu$ and for all $r < R^{-1}$, 
then
\begin{equation} \label{secondestimate}\int_0^{\pi} \int_{\mathbb{R} } \left\lvert \pi_{\sharp } P_{\mathbb{V}_{\theta \sharp}^{\perp} }\left(\mu \ast \widecheck{\psi_R} \right) \right\rvert^p \, d\mathcal{H}^1_E \, d\theta \leq C_{\epsilon}  C^{p-1}  \mu(\mathbb{H}) R^{\epsilon -(\alpha-1)(p-1)}. \end{equation}\end{theorem}
\begin{remark} To get $p=4/3$ requires the local smoothing inequality from \cite{gaoliumiaoxi}, but the local smoothing inequality from \cite{beltranhickmansogge} would be sufficient for $1 < p \leq 6/5$, and any $p>1$ would suffice for the applications to projections below. \end{remark}
\begin{remark} The simpler inequality \eqref{secondestimate} will be used for the packing dimension bound, but the more precise \eqref{thirdestimate} will be needed for the Hausdorff dimension bound. \end{remark}
\begin{sloppypar}\begin{proof} The inequality \eqref{firstestimate} will be proved first, and then the minor changes to the proof of \eqref{firstestimate} necessary for \eqref{thirdestimate} (which implies \eqref{secondestimate}) will be explained. 

By approximation (using that the dual of $L^p$ has a dense subset of $C_0^{\infty}$ functions when $p>1$), it suffices to prove \eqref{firstestimate} under the assumption that $\mu \in C_0^{\infty}(\mathbb{R}^3)$.

Let
\[ Af(z,t) = \int_0^{\pi} \left( \chi f \right)\left( \theta, t + \frac{\omega\left( \pi_{V_{\theta}}(z), z \right) }{2} \right) \, d\theta, \]
where $\chi$ is a smooth bump equal to 1 on $[0,\pi] \times J$ and vanishing on a slightly larger rectangle, where $J$ is an interval of length $\sim 1$. 

By Lemma~\ref{adjoint} and duality, it suffices to prove that for any smooth compactly supported function $f$,
\begin{equation} \label{needtobound}  \left\lvert \int_{\mathbb{H}}  Af(z,t) \, d\mu(z,t) \right\rvert \leq C_{\epsilon, \alpha}  c_{\alpha,E}(\mu)^{1/p'} \mu(\mathbb{H})^{1/p} R^{\epsilon-(\alpha-1)/p'} \lVert f \rVert_{p'}, \end{equation}
where $p'$ is the Hölder conjugate of $p$.  Fix such an $f$ and decompose 
\begin{equation} \label{littlewoodpaley} f = f_0 + \sum_{0 < k < \log_2 R} f_k + \sum_{k \geq \log_2 R}  f_k, \end{equation}
where $f_k$ is frequency supported in $|\xi| \sim 2^{k}$ for $k \geq 1$, and $f_k = f \ast \widecheck{\phi_k}$ with $\phi_k$ a smooth bump on $|\xi| \sim 2^k$. The term $f_{0}$ is $f_{0} = f \ast \widecheck{\psi}$, with $\psi$ a smooth bump on $|\xi| \lesssim 1$. If the term from $f_0$ dominates the left-hand side of \eqref{needtobound}, then 
\[ \lVert Af_{0} \rVert_{\infty} \lesssim \lVert f_{0} \rVert_{\infty} \lesssim \lVert f \rVert_{p'}, \]
and thus, since $\mu$ is supported in a Euclidean ball of radius $R^{-1}$,
\[ \left\lvert \int_{\mathbb{H}}  Af_{0}(z,t) \, d\mu(z,t) \right\rvert \lesssim \mu(\mathbb{H}) \lVert f \rVert_{p'}\leq \mu(\mathbb{H})^{1/p} c_{\alpha,E}(\mu)^{1/p'}R^{-\alpha/p'}\lVert f \rVert_{p'}, \]
which is better than \eqref{needtobound}. 

For the remaining frequencies, by summing two geometric series, 
it suffices to show that for any positive integer $k$ and sufficiently small $\epsilon >0$,
\begin{multline} \label{aim} \left\lvert \int_{\mathbb{H}} \mu(z,t) Af_k(z,t) \, dz \, dt\right\rvert 
 \\
\leq C_{\epsilon} 2^{k\epsilon} 2^{k/p'} \mu(\mathbb{H})^{1/p}  c_{\alpha,E}(\mu)^{1/p'} \min\left\{  2^{-k\alpha/p'}, R^{-\alpha/p'} \right\} \lVert f \rVert_{p'}.
\end{multline}
Let $k$ be given. Let $B$ be a Euclidean ball of radius $\sim 1$ containing the support of $\mu$, with $|z| \sim 1$ for all $(z,t) \in B$. 

For each $t \in \mathbb{R}$, define $\Phi: \mathbb{R}^2 \times \mathbb{R}^2 \to \mathbb{R}$ by 
\[ \Phi_t(z, \theta, r) = t + \frac{1}{2} \omega\left(\pi_{V_{\theta}}(z), z \right) - r. \]
By writing $z = x_1+ix_2$ and using the definition in \cite[p.~494]{stein}, the rotational curvature of $\Phi_t$ is 
\[ \rot \curv \Phi_t = \det \begin{pmatrix} \Phi_t & \partial_{\theta} \Phi_t & -1 \\ \partial_{x_1} \Phi_t & \partial_{\theta x_1} \Phi_t & 0 \\ \partial_{x_2} \Phi_t & \partial_{\theta x_2} \Phi_t & 0 \end{pmatrix}. \]
A formula for $\Phi_t$ is 
\[ \Phi_t(x_1,x_2, \theta, r )  = t + \frac{1}{2} \left( x_1 \cos \theta + x_2 \sin \theta \right) \left( x_2 \cos \theta - x_1 \sin \theta \right) - r. \]
Hence 
\[ \partial_{x_1} \Phi_t = \frac{1}{2} \left( x_2 \cos(2\theta) - x_1 \sin(2\theta) \right), \]
and 
\[ \partial_{x_2} \Phi_t = \frac{1}{2} \left( x_1 \cos(2\theta) + x_2 \sin(2\theta) \right). \]
This gives 
\begin{equation} \label{identity1} \partial_{\theta x_1} \Phi_t = -2 \partial_{x_2} \Phi_t, \end{equation}
and 
\begin{equation} \label{identity2} \partial_{\theta x_2} \Phi_t = 2 \partial_{x_1} \Phi_t. \end{equation}
Hence 
\[ \rot \curv \Phi_t  = -2 \left[ \left( \partial_{x_1} \Phi_t \right)^2 + \left( \partial_{x_2} \Phi_t \right)^2 \right] = -(x_1^2+x_2^2)/2. \]
Therefore $\left\lvert \rot \curv \Phi_t(z,\theta, r) \right\rvert \sim 1$ for $(z,t)\in B$. It follows from \cite[p.~496 and \S~4.8(a) on p.~517]{stein} that for each fixed $t \in \mathbb{R}$, $f \mapsto Af(\cdot, t )$ is a Fourier integral operator of order $-1/2$. 

To verify the cinematic curvature condition from \cite{sogge}, by the above, either $\left\lvert \partial_{x_1} \Phi_t \right\rvert \sim 1$ or $\left\lvert \partial_{x_2} \Phi_t \right\rvert \sim 1$ for $(z,t) \in B$. By rotation invariance, it may be assumed that $\left\lvert \partial_{x_2} \Phi_t \right\rvert \sim 1$. Then by \cite[Theorem~2.1]{kung}, the ``cinematic curvature'' of the operator $f \mapsto Af$ (defined as $\cin \curv$ in \cite{kung}) is (for $(z,t) \in B$)
\begin{equation} \label{cincurvformula} \cin \curv  \sim  \det\begin{pmatrix} \partial_{x_1} \Phi_t & \partial_{x_2} \Phi_t & 1 \\ \partial_{\theta x_1} \Phi_t & \partial_{\theta x_2} \Phi_t & 0 \\ \partial_{\theta \theta x_1} \Phi_t & \partial_{\theta \theta x_2} \Phi_t & 0 \end{pmatrix}. \end{equation}
More precisely, Theorem~2.1 from \cite{kung} is that the cinematic curvature condition from \cite{sogge} for the operator $f \mapsto Af$ is equivalent to the nonvanishing of the quantity $\cin \curv$ defined above, for $(z,t) \in B$. By \eqref{identity1}, \eqref{identity2}, and \eqref{cincurvformula}, 
\begin{multline*} \cin \curv \sim \partial_{\theta x_1} \Phi_t \partial_{\theta \theta x_2} \Phi_t - \partial_{\theta x_2 } \Phi_t \partial_{\theta \theta x_1} \Phi_t  \\
= 4 \left( \left( \partial_{x_1} \Phi_t \right)^2 + \left( \partial_{x_2} \Phi_t \right)^2 \right) = x_1^2+x_2^2 \sim 1, \end{multline*}
for $(z,t) \in B$. This verifies the cinematic curvature condition for the operator $f \mapsto Af$ in $B$, and that the operators $Af(\cdot, t)$ are Fourier integral operators of order $-1/2$.  Therefore, by the variable coefficient local smoothing inequality (\cite[Theorem~1.4 with $\mu = -1/2$]{gaoliumiaoxi} for $p' \geq 4$ or alternatively \cite{beltranhickmansogge} for $p' \geq 6$),  for any $\epsilon >0$,
\begin{equation} \label{eq:gaoliumiaoxi} \lVert Af_k \rVert_{L^{p'}(B)} \leq C_{\epsilon} 2^{k \epsilon} 2^{-2k/p'} \lVert f \rVert_{p'}, \end{equation}
For $N \in \mathbb{N}$ and $\partial \in \{\partial_{x_1}, \partial_{x_2}, \partial_t \}$, $\partial^N Af_k$ is a weighted sum of similar averaging operators to $A$ applied to derivatives of $f_k$ up to order $N$. Therefore, similarly to \eqref{eq:gaoliumiaoxi}, for any $\epsilon >0$,
\begin{equation} \label{derivativeinequality} \lVert  \partial^N Af_k  \rVert_{L^{p'}(B)} \leq C_{N, \epsilon} 2^{k \epsilon} 2^{-2k/p'} \lVert f_k \rVert_{W^{N,p'}}. \end{equation}
The gain of $2^{-2k/p'}$ in \eqref{derivativeinequality} will not be needed, so the local smoothing inequality \eqref{derivativeinequality} could be replaced by an interpolation of the simpler $L^2$ and $L^{\infty}$ bounds; it is just used here to simplify the referencing. By Young's convolution inequality,
\[ \lVert f_k \rVert_{W^{N,p'}} \leq C_N \lVert f \rVert_{p'} \left\lVert \widecheck{\phi_k} \right\rVert_{W^{N,1}} \leq C_N \lVert f \rVert_{p'} 2^{kN}. \] 
Hence 
\[ \lVert  \partial^N Af_k  \rVert_{L^{p'}(B)} \leq C_N 2^{kN} \lVert f \rVert_{p'}, \]
where the factor $2^{k \epsilon} 2^{-2k/p'}$ has been removed as it provides no benefit here. Integrating by parts many times and applying Hölder's inequality yields that $\widehat{\chi_B Af_k}$ is rapidly decaying outside $B_E(0, 2^k)$, where $\chi_B$ is a smooth bump function adapted to $B$. Hence
\begin{multline} \label{afterholder} \left\lvert \int_{\mathbb{H}} \mu(z,t) Af_k(z,t) \, dz \, dt\right\rvert \leq  \int_{B} |\mu \ast \psi_k(z,t)| |Af_k(z,t)| \, dz \, dt \\
+ C_{\epsilon} 2^{-100k} \lVert f \rVert_{p'} \mu(\mathbb{H}), \end{multline}
where $\psi_k$ is a non-negative smooth bump function, with $\psi_k \sim 2^{3k(1+\epsilon)}$ on $B_3(0, 2^{-k(1+\epsilon)})$ and rapidly decaying outside this ball, with $\int \psi_k \lesssim 1$. By substituting into \eqref{aim}, it remains to show that 
\begin{multline} \label{aim2} \int_{B} \left\lvert \mu \ast \psi_k(z,t)\right\rvert |Af_k(z,t)| \, dz \, dt 
 \\
\leq C_{\epsilon} 2^{k\epsilon} 2^{k/p'} \mu(\mathbb{H})^{1/p} c_{\alpha,E}(\mu)^{1/p'} \min\left\{ 2^{-k\alpha/p'}, R^{-\alpha/p'}\right\} \|f\|_{p'}.  
\end{multline}
By Hölder's inequality,
\[  \int_{B} \left\lvert \mu \ast \psi_k (z,t)\right\rvert \left\lvert Af_k(z,t) \right\rvert\, dz \, dt \leq \lVert \mu \ast \psi_k \rVert_p \lVert Af_k \rVert_{L^{p'}(B)}. \]
Applying \eqref{eq:gaoliumiaoxi} to the above gives, for any $\epsilon >0$, 
\begin{equation} \label{aftersmoothing} \int_{B} \left\lvert \mu \ast \psi_k (z,t)  \right\rvert \left\lvert Af_k(z,t) \right\rvert\, dz \, dt \lesssim \left\lVert \mu \ast \psi_k \right\rVert_p 2^{-k\left(\frac{2}{p'} - \epsilon\right)} \left\lVert f \right\rVert_{p'}. \end{equation}
Since $\mu$ is supported in a Euclidean ball of radius $R^{-1}$, and since $\psi_k$ is rapidly decaying outside $B_E(0, 2^{-k})$, the Euclidean Frostman condition on $\mu$ gives
\begin{equation} \label{linfinity} \left\lVert \mu \ast \psi_k \right\rVert_{\infty} \lesssim 2^{k(3+O(\epsilon))} \min\left\{  2^{-k\alpha}, R^{-\alpha}  \right\} c_{\alpha,E}(\mu). \end{equation}
Hence
\begin{multline*} \left\lVert \mu \ast \psi_k \right\rVert_p^{p} \lesssim \lVert \mu \ast \psi_k \rVert_{\infty}^{p-1} \mu(\mathbb{H}) \\
\lesssim 2^{k(p-1)(3 + O(\epsilon))} \min\left\{ 2^{-k\alpha(p-1)} R^{-\alpha(p-1)} \right\} c_{\alpha,E}(\mu)^{p-1} \mu(\mathbb{H}). \end{multline*}
Substituting into \eqref{aftersmoothing} gives 
\begin{multline*} \int_{B} \left\lvert \mu \ast \psi_k (z,t)  \right\rvert \left\lvert Af_k(z,t) \right\rvert\, dz \, dt \lesssim \\
\left(2^{k(p-1)(3+ O(\epsilon))} \min\left\{ 2^{-k\alpha(p-1)} R^{-\alpha(p-1)} \right\} c_{\alpha,E}(\mu)^{p-1} \mu(\mathbb{H})\right)^{1/p} 2^{-k\left(\frac{2}{p'} - \epsilon\right)}   \lVert f \rVert_{p'}\\
= 2^{k\left( \frac{1}{p'} + O(\epsilon) \right)} \min\left\{ 2^{-k\alpha/p'}, R^{-\alpha/p'} \right\} c_{\alpha,E}(\mu)^{1/p'} \mu(\mathbb{H})^{1/p} \lVert f \rVert_{p'}. \end{multline*} 
This verifies \eqref{aim2} and finishes the proof of \eqref{firstestimate}. 

For the proof of \eqref{thirdestimate} (which implies \eqref{secondestimate}), the main difference is that instead of \eqref{needtobound} it is required to show that
\begin{multline} \label{alteredversion} \left\lvert \int Af(z,t) \left( \mu \ast \widecheck{\psi_R} \right)(z,t) \, dz \, dt \right\rvert \\ \leq C_{\epsilon} \mu(\mathbb{H})^{1/p} \sum_{k \geq \log_2 R } 2^{k \epsilon} \left( \sup_{x \in \mathbb{H}} 2^k \mu\left(B_E(x,2^{-k}\right) \right)^{1/p'} \lVert f \rVert_{p'}. \end{multline}
 Since, as explained previously, $\widehat{\chi_B Af_k}$ is rapidly decaying outside $B_E(0, 2^k)$, where $\chi_B$ is a smooth bump function on a Euclidean ball $B$ of radius $\sim 1$ containing the support of $\mu$ and with $|z| \sim 1$ for all $(z,t) \in B$, the only frequencies in the decomposition \eqref{littlewoodpaley} contributing non-negligibly to the left-hand side of \eqref{alteredversion} are those with $2^{k} \geq R^{1-\epsilon}$. Therefore, it suffices to show that for $k \geq (1-\epsilon)\log_2 R$, 
\begin{multline} \label{newaim} \int_B \left\lvert Af_k(z,t) \right\rvert \left\lvert \left(\mu \ast \psi_k \right)(z,t)\right\rvert \, dz \, dt  \\
\leq C_{\epsilon}  \mu(\mathbb{H})^{1/p} 2^{kO(\epsilon)} \left( \sup_{x \in \mathbb{H}} 2^k \mu\left(B_E(x,2^{-k}\right) \right)^{1/p'} \lVert f \rVert_{p'}, \end{multline}
where $B$ is a unit ball with $|z| \sim 1$ for all $(z,t) \in B$, and $\psi_k$ is as in \eqref{afterholder}. The reason that the range $(1-\epsilon) \log_2 R \leq k < \log_2 R$ does not need to be included in \eqref{alteredversion} is that this range is absorbed into the $\log_2 R \leq k < \log_2 R +1$ term by the $2^{k \epsilon}$ factor occurring in \eqref{alteredversion}. Due to the restriction $2^{k} \geq R^{1-\epsilon}$, the proof of \eqref{newaim} is nearly identical to the proof of \eqref{aim2} in this case, since the only frequencies which made significant use of the support of $\mu$ having Euclidean diameter $\lesssim R^{-1}$ in \eqref{linfinity} were for $2^k \leq R^{1-\epsilon}$. The only modification to the proof necessary is that instead of using the Euclidean Frostman condition in \eqref{linfinity}, the bound 
\[ \|\mu \ast \psi_k\|_{\infty} \lesssim 2^{kO(\epsilon)} 2^{3k} \sup_{x \in \mathbb{H}} \mu(B_E(x,2^{-k} ) ) \]
should be used instead. Apart from this change,  following the rest of the proof of \eqref{aim2} results in \eqref{newaim}.    \end{proof}\end{sloppypar}

\section{An intersection theorem} 

Recall that $\pi: \mathbb{H} \to \mathbb{R}$ is $\pi(z,t) = t$. The Hausdorff dimension version of the lemma below is the planar case of \cite[Lemma~3.2]{mattila2}, though in \cite{mattila2} the author states that the planar case is essentially due to Marstrand \cite[Lemma~16]{marstrand}. Below, $\mathcal{P}^t_E$ refers to the Euclidean $t$-dimensional packing measure. 

\begin{lemma} \label{analysis} Fix $\theta \in [0, \pi)$. Let $F \subseteq \mathbb{V}_{\theta}^{\perp}$ be a Borel set, and $t>0$. 
\begin{enumerate}
\item If $\mathcal{H}^t_E\left( F \cap \pi^{-1}(u) \right) = 0$ for all $u \in \mathbb{R}$, then for any finite Borel measure $\nu$ on $\mathbb{V}_{\theta}^{\perp}$, 
\[ \limsup_{r \to 0^+} \liminf_{\delta \to 0^+} r^{-t} \delta^{-1} \nu\left\{ y \in B_E(x,r) : d_{E}\left(\pi(x), \pi(y) \right) < \delta \right\} = \infty, \]
for $\nu$-a.e.~$x \in F$.
\item  If $\mathcal{P}^t_E\left( F \cap \pi^{-1}(u) \right) =0$ for all $u \in \mathbb{R}$, then for any finite Borel measure $\nu$ on $\mathbb{V}_{\theta}^{\perp}$, 
\[ \liminf_{r \to 0^+} \liminf_{\delta \to 0^+} r^{-t} \delta^{-1} \nu\left\{ y \in B_E(x,r) : d_{E}\left(\pi(x), \pi(y) \right) < \delta \right\} = \infty, \]
for $\nu$-a.e.~$x \in F$.   \end{enumerate}  \end{lemma} 

\begin{proof} Since the Hausdorff measure version was proved in \cite{mattila2}, only the packing measure version will be proved here. The proof is similar to that in \cite{mattila2}. Since any finite Borel measure on Euclidean space is inner regular, it suffices to show that $\nu(E) = 0$ for any positive integer $N$, for any compact set $E$ with 
\begin{multline*} E \subseteq \\
\left\{ x \in F : \liminf_{r \to 0^+} \liminf_{\delta \to 0^+} r^{-t} \delta^{-1} \nu\left\{ y \in B_E(x,r) : d_{E}\left(\pi(x), \pi(y) \right) < \delta \right\} \leq N \right\}.  \end{multline*} 
Let $\mu$ be the restriction of $\nu$ to $E$, given by $\mu(A) = \mu(A \cap E)$, so that $\supp \mu \subseteq E$ (as $E$ is compact). 

It will be shown that $\pi_{\sharp}\mu \ll \mathcal{H}^1_E$. For this, it suffices to show that $\pi_{\sharp}\left( \nu \chi_{G} \right) \ll \mathcal{H}^1_E$, whenever $r>0$ is fixed and $G$ is a compact subset of 
\[ \left\{ x \in F : \liminf_{\delta \to 0^+} \delta^{-1} \nu\left\{ y \in B_E(x,r) : d_{E}\left(\pi(x), \pi(y) \right) < \delta \right\} \leq 2N \right\}. \]
Cover $G$ by boundedly overlapping closed balls $B$ of radius $r/2$ and let $G_B = G \cap B$. It suffices to prove that $\pi_{\sharp}\left( \nu \chi_{G_B} \right) \ll \mathcal{H}^1_E$. If $\lambda \in \supp\left( \pi_{\sharp} \left( \nu \chi_{G_B} \right)\right)$, let $x_0 \in G_B$ be such that $\pi(x_0) = \lambda$. Then, by the definition of $G$,
\begin{align*} \liminf_{\delta \to 0^+} \delta^{-1} \pi_{\sharp} \left( \nu \chi_{G_B} \right)\left(B(\lambda, \delta )\right) &= \liminf_{\delta \to 0^+} \delta^{-1} \nu\left\{ y \in G_B : d_E(\pi(x_0), \pi(y) ) < \delta \right\} \\
&\leq \liminf_{\delta \to 0^+} \delta^{-1} \nu\left\{ y \in B(x_0,r) : d_E(\pi(x_0), \pi(y) ) < \delta \right\} \\
& \leq 2N. \end{align*}
Therefore, for any $\lambda \in \supp\left( \pi_{\sharp} \left( \nu \chi_{G_B} \right)\right)$, there exists arbitrarily small $\delta>0$ such that $\pi_{\sharp} \left( \nu \chi_{G_B} \right)\left(B(\lambda, \delta )\right) \leq 3N \delta$. By the Vitali covering lemma and since Lebesgue measure is outer regular, it follows that $\pi_{\sharp} \left( \nu \chi_{G_B} \right) \ll \mathcal{H}^1_E$, and as explained above this yields that $\pi_{\sharp}  \mu \ll \mathcal{H}^1_E$. 

Suppose for a contradiction that $\mu(E) =\nu(E) >0$. By \cite[Chapter~10]{mattila}, for $\mathcal{H}^1_E$-a.e.~$\lambda \in \mathbb{R}$, there exists a finite Borel measure $\mu_{\lambda}$ with
\begin{equation} \label{supportcondition} \supp \mu_{\lambda} \subseteq \supp \mu \cap \pi^{-1}(\lambda) \subseteq E \cap \pi^{-1}(\lambda), \end{equation}
 such that 
\begin{equation} \label{slicemeasure} \int \phi \, d\mu_{\lambda} = \lim_{\delta \to 0^+} (2\delta)^{-1} \int_{\left\{ (z,y) \in \mathbb{V}_{\theta}^{\perp} : |y-\lambda| \leq \delta\right\}} \phi \, d\mu, \end{equation}
for all non-negative continuous functions $\phi$ on $\mathbb{V}_{\theta}^{\perp}$. Therefore, by Chebychev's inequality, there exists $\lambda$ such that \eqref{supportcondition} and \eqref{slicemeasure} hold, and such that
\[  \lim_{\delta \to 0^+} (2\delta)^{-1} (\pi_{\sharp}\mu)(B_E(\lambda, \delta) ) = (\pi_{\sharp}\mu)(\lambda) \gtrsim C^{-1} \mu(E), \]
where $C = \mathcal{H}^1_E(\pi(E))$, which is nonzero by the assumption $\mu(E)>0$ and since $\pi_{\sharp}\mu \ll \mathcal{H}^1_E$. By taking $\phi=1$ in \eqref{slicemeasure}, this yields 
\[ \mu_{\lambda}(\mathbb{R}) \gtrsim C^{-1} \mu(E), \]
and in particular $\mu_{\lambda}(\mathbb{R}) >0$. For any $(z_0,\lambda) \in \supp \mu_{\lambda}$, 
\begin{multline*} \liminf_{r \to 0^+} r^{-t} \mu_{\lambda}(B_E((z_0,\lambda),r) ) \lesssim \liminf_{r \to 0^+} r^{-t} \mu_{\lambda}(B_E((z_0,\lambda),r/10) )\\
\leq \liminf_{r \to 0^+} r^{-t} \int \phi_r \, d\mu_{\lambda}, \end{multline*}
where $\phi_r$ is a smooth bump function with $0 \leq \phi_r \leq 1$, which is equal to $1$ on $B_E((z_0,\lambda),r/10)$ and which vanishes outside $B_E((z_0,\lambda),r)$. By \eqref{slicemeasure}, it follows that
\begin{multline*} \liminf_{r \to 0^+} r^{-t} \mu_{\lambda}(B_E((z_0,\lambda),r) ) \\
\lesssim \liminf_{r \to 0^+} \liminf_{\delta \to 0^+} \delta^{-1} r^{-t} \mu\left\{ x \in B_E((z_0,\lambda),r) : |\pi(x)-\pi(z_0,\lambda)| < \delta \right\}.  \end{multline*}
But $(z_0,\lambda) \in E$ by \eqref{supportcondition}, so this yields that 
\[ \liminf_{r \to 0^+} r^{-t} \mu_{\lambda}(B_E((z_0,\lambda),r) ) \lesssim N. \]
By the mass distribution principle for packing dimension (see e.g.~\cite[Chapter~2, Proposition~2.2(c)]{falconertechniques}
), it follows that $\mathcal{P}^t_E(E \cap \pi^{-1}(\lambda)) >0$, which contradicts the assumption.
\end{proof}

The Korányi metric equals the Euclidean metric on the intersection of any fibre of $\pi$ with a vertical plane (any line of constant height inside a vertical plane), so the Euclidean Hausdorff measure $\mathcal{H}^t_E$ in Lemma~\ref{analysis} could be replaced by the Korányi Hausdorff measure, and similarly for the packing measure. 

The following theorem consists of the two key intersection results which will imply Theorem~\ref{etccorollary1} and Proposition~\ref{etccorollary}.

\begin{theorem}  \label{intersectiontheorem2} Let $2 < s < 3$. 
\begin{enumerate} 
\item Suppose that $\mu$ is a compactly supported Borel measure on $\mathbb{H}$ which is Euclidean Ahlfors-regular. Then for any Borel set $A \subseteq \supp \mu$ with $\mu(A)>0$ such that $c_{s,\mathbb{H}}(\mu \restriction_A)<\infty$, for a.e.~$\theta \in [0, \pi)$, 
\begin{equation} \label{firstclaim} \mathcal{H}^1_E\left\{ \lambda \in \mathbb{R} : \dim\left(\pi^{-1}(\lambda) \cap P_{\mathbb{V}_{\theta}^{\perp}}(A)\right) \geq s-2 \right\} >0. \end{equation}
\item Suppose that $A \subseteq \mathbb{H}$ is a Borel (or analytic) with $\dim A = s$. Then for any $\epsilon >0$, for a.e.~$\theta \in [0, \pi)$, 
\begin{equation} \label{secondclaim} \mathcal{H}^1_E\left\{ \lambda \in \mathbb{R} : \dim_P\left(\pi^{-1}(\lambda) \cap P_{\mathbb{V}_{\theta}^{\perp}}(A)\right) \geq s-2-\epsilon \right\} >0. \end{equation} \end{enumerate}  \end{theorem} 
\begin{proof} The proof of the second statement \eqref{secondclaim} has some additional technical steps compared to the first statement \eqref{firstclaim}, so only the proof of the second statement will be given in full, and then the minor adjustments and simplifications required to prove the first part will be explained. 

 By Heisenberg dilation, vertical translation, and since $s> 2$, it may be assumed that $A$ is contained in a set of the form
\[ \{ (z,t) : 1 \leq |z| \leq 2, |t| \leq 1 \}. \]
Fix such a set $A$. Since the conclusion \eqref{secondclaim} allows $\epsilon$ losses, by Frostman's lemma (\cite{howroyd}) it may be assumed that $A$ supports a finite Borel measure with $c_{s, \mathbb{H}}(\mu) < \infty$. By scaling it may be assumed that $c_{s, \mathbb{H}}(\mu)=1$, which will simplify the notation in a couple of places. 

 Let $ 0 < t < s-2$. The projection inequalities from Lemma~\ref{projectionconjecture} and Theorem~\ref{lpestimate} will be used to show that, for some $p>1$ possibly depending on $t$ and $s$, and for some Borel set $F \subseteq \supp \mu$ with $\mu(F) >0$ depending on $t$,
\begin{multline} \label{desiredinequality} \int_{0}^{\pi} \int \liminf_{r \to 0^+} \liminf_{\delta \to 0^+} \left(r^{-t} \delta^{-1} \left(P_{\mathbb{V}_{\theta \sharp}^{\perp}}\mu \right)\left\{ y \in B_E(x,r) : d_E(\pi(x), \pi(y) )< \delta \right\}\right)^{p-1}  \\
 d\left( P_{\mathbb{V}_{\theta \sharp}^{\perp}}\mu_F \right)(x) \, d\theta =0; \end{multline}
the value of $p$ not being important for the application to intersections below. It will first be shown that \eqref{desiredinequality} implies \eqref{secondclaim} in the theorem. Assuming \eqref{desiredinequality}, for a.e.~$\theta \in [0, \pi)$,
\begin{multline} \label{desired2} \int \liminf_{r \to 0^+} \liminf_{\delta \to 0^+} \\
\left(r^{-t} \delta^{-1} \left(P_{\mathbb{V}_{\theta \sharp}^{\perp}}\mu \right)\left\{ y \in B_E(x,r) : d_E(\pi(x), \pi(y) )< \delta \right\}\right)^{p-1} \, d\left( P_{\mathbb{V}_{\theta \sharp}^{\perp}}\mu_F \right)(x) = 0, \end{multline}
and $\pi_{\sharp} P_{\mathbb{V}_{\theta \sharp}^{\perp}} \mu \ll \mathcal{H}^1_E$ (by Theorem~\ref{lpestimate}, using $s>2$ and dimension comparison \eqref{comparisonfrostman} below). For such a $\theta$, let 
\[ G_{\theta} = \left\{ \lambda \in \mathbb{R} : \mathcal{P}^t_E\left(\pi^{-1}(\lambda) \cap P_{\mathbb{V}_{\theta}^{\perp}}(A)\right) =0 \right\}. \]
By defining $B = P_{\mathbb{V}_{\theta}^{\perp}}(A) \cap \pi^{-1}(G_{\theta})$, it is straightforward to check that $\mathcal{P}^t_E(\pi^{-1}(\lambda) \cap B) = 0$ for every $\lambda \in \mathbb{R}$. Hence, by Lemma~\ref{analysis} and since $P_{\mathbb{V}_{\theta \sharp}^{\perp}} \mu$ has support contained in $P_{\mathbb{V}_{\theta}^{\perp}}(A)$, it holds that for $P_{\mathbb{V}_{\theta \sharp}^{\perp}} \mu$-a.e.~$x \in \pi^{-1}(G_{\theta})$,
\[ \liminf_{r \to 0^+} \liminf_{\delta \to 0^+} r^{-t} \delta^{-1} \left(P_{\mathbb{V}_{\theta \sharp}^{\perp}}\mu\right)\left\{ y \in B_E(x,r) : d_{E}\left(\pi(x), \pi(y) \right) < \delta \right\} = \infty. \]
Comparing with \eqref{desired2} gives that
\[ \left(\pi_{\sharp} P_{\mathbb{V}_{\theta \sharp}^{\perp}}\mu_F\right)(G_{\theta}) = \left(P_{\mathbb{V}_{\theta \sharp}^{\perp}}\mu_F\right)(\pi^{-1}(G_{\theta})) = 0, \]
for a.e.~$\theta \in [0, \pi)$. It follows that for a.e.~$\theta \in [0, \pi)$, $\mathcal{P}^t_E\left( \pi^{-1}(\lambda) \cap P_{\mathbb{V}_{\theta}^{\perp}}(A) \right) >0$ for $\pi_{\sharp} P_{\mathbb{V}_{\theta \sharp}^{\perp}}\mu_F$-a.e.~$\lambda \in \mathbb{R}$. Since this holds for any $t < s-2$, and since $\mu(F) >0$ and (for a.e.~$\theta$) $\pi_{\sharp} P_{\mathbb{V}_{\theta \sharp}^{\perp}} \mu \ll \mathcal{H}^1_E$, it implies that for any $t < s-2$, for a.e.~$\theta \in [0, \pi)$, 
\[ \mathcal{H}^1_E\left\{ \lambda \in \mathbb{R} : \dim_P\left( \pi^{-1}(\lambda) \cap P_{\mathbb{V}_{\theta}^{\perp}}(A) \right) \geq t \right\} >0, \]
as claimed.

It remains to prove \eqref{desiredinequality}, for any $0 < t < s-2$. Let $\epsilon = \frac{1}{10^4}\left[s-2-t\right]>0$. By pigeonholing, there exists $s_E \in (1, 3]$ and a Borel set $F \subseteq \supp \mu$ with $\mu(F) >0$ such that for all $x \in F$, 
\begin{equation} \label{lowerdensity} s_E - \epsilon < \liminf_{r \to 0} \frac{ \log \mu\big(B_E(x,r)\big) }{\log r} < s_E + \epsilon. \end{equation}
The assumption $s>2$ together with the dimension comparison inequality \eqref{comparisonfrostman} below ensures that $s_E>1$ (provided $\epsilon$ is sufficiently small, or equivalently if $t$ is sufficiently close to $s-2$). By \eqref{lowerdensity}, for all $x \in F$, 
\[  \mu\big(B_E(x,r)\big) < r^{s_E-\epsilon} \qquad \text{ for all } r < r_{0}, \]
where $r_{0}>0$ is independent of $x$ (after replacing $F$ by a Borel subset of positive measure). Similarly, by \eqref{lowerdensity}, for all $x \in F$, there exists arbitrarily small $r >0$ with 
\[ \mu\big(B_E(x,r)\big) > r^{s_E+\epsilon}. \]
To prove \eqref{desiredinequality}, by the monotone convergence theorem, it suffices to show that
\begin{multline} \label{monotoneconv} \lim_{K \to \infty}  \int_{0}^{\pi} \int \inf_{k \geq K} \liminf_{\delta \to 0^+} \\
\left(2^{kt} \delta^{-1} \left(P_{\mathbb{V}_{\theta \sharp}^{\perp}}\mu \right)\left\{ y \in B_E\left(x,2^{-k}\right) : d_E(\pi(x), \pi(y) )< \delta \right\}\right)^{p-1} \\
d\left( P_{\mathbb{V}_{\theta \sharp}^{\perp}}\mu_F \right)(x) \, d\theta =0. \end{multline}
Fix a positive integer $K$ with $2^{-K} < r_{0}$. For each $x \in F$, write $k=k(x)$ for the smallest positive integer with $k \geq K$ such that $\mu\big(B_E(x,2^{-k})\big) > 2^{-k(s_E+\epsilon)}/100$, which exists by \eqref{lowerdensity}. Decompose $F = \bigcup_{k \geq K} F_k$ according to value of $k=k(x)$. To prove \eqref{monotoneconv}, it suffices to show that 
\begin{multline}   \label{sumbound}\sum_{k \geq K} \int_{0}^{\pi} \int \liminf_{\delta \to 0^+} \\
\left(2^{kt} \delta^{-1} \left(P_{\mathbb{V}_{\theta \sharp}^{\perp}}\mu \right)\left\{ y \in B_E\left(x,2^{-k}\right) : d_E(\pi(x), \pi(y) )< \delta \right\}\right)^{p-1}  \\
 d\left( P_{\mathbb{V}_{\theta \sharp}^{\perp}}\mu_{F_k} \right)(x) \, d\theta \lesssim 2^{-K\epsilon (p-1)} \mu(\mathbb{H}) c_{s,\mathbb{H}}(\mu)^{p-1}. \end{multline}
For this, it is enough to show that for any $k \geq K$, 
\begin{multline*}  \int_{0}^{\pi} \int \liminf_{\delta \to 0^+} \left(2^{kt} \delta^{-1} \left(P_{\mathbb{V}_{\theta \sharp}^{\perp}}\mu \right)\left\{ y \in B_E\left(x,2^{-k}\right) : d_E(\pi(x), \pi(y) )< \delta \right\}\right)^{p-1}  \\
 d\left( P_{\mathbb{V}_{\theta \sharp}^{\perp}}\mu_{F_k} \right)(x) \, d\theta \lesssim 2^{-k\epsilon (p-1)} \mu(\mathbb{H}) c_{s,\mathbb{H}}(\mu)^{p-1}. \end{multline*}
By Fatou's lemma, it suffices to find, for any $\epsilon >0$, a $p>1$ depending only on $s$ and $\epsilon$, such that for any positive integer $k$ and any $\delta>0$,
 \begin{multline} \label{reusedbelow11} \int_0^{\pi} \int 
 \left(\delta^{-1} \left(P_{\mathbb{V}_{\theta \sharp}^{\perp}}\mu \right)\left\{ y \in B_E\left(x,2^{-k}\right) : d_E(\pi(x), \pi(y) )< \delta \right\}\right)^{p-1}  \\
 d\left( P_{\mathbb{V}_{\theta \sharp}^{\perp}}\mu_{F_k} \right)(x) \, d\theta \lesssim \mu(\mathbb{H}) c_{s,\mathbb{H}}(\mu)^{p-1} 2^{-k(p-1) \left( s-2 - O(\epsilon)\right)}, \end{multline}
for any Borel measure $\mu$ with $c_{s,\mathbb{H}}(\mu)=1$, supported in
\[ \{ (z,t) \in \mathbb{H} : 1 \leq |z| \leq 2, |t| \leq 1\}, \]
and whenever $F_k \subseteq \supp \mu$ satisfies 
\begin{equation} \label{lowerdensity2} \mu\big(B_E(x,2^{-k})\big) \gtrsim 2^{-k(s_E+\epsilon )}, \qquad x \in \supp F_k, \end{equation}
and 
\begin{equation} \label{upperdensity2} \mu\big(B_E(x, r)\big) \leq r^{s_E-\epsilon}, \qquad x \in \supp F_k, \quad r < 2^{-k}. \end{equation}
In \eqref{reusedbelow11}, $O(\epsilon)$ can be taken as $1000\epsilon$. Let $\eta>0$ be very small, to be chosen after $\epsilon$ but before $p$, and assume that $\epsilon >0$ is very small.

Let $\mathcal{B}$ be a boundedly overlapping cover of the support of $F_k$ by Euclidean balls of radius $2^{-k}$. Then
 \begin{multline} \label{sumoverB} \int_0^{\pi} \int 
 \left(\delta^{-1}\left(P_{\mathbb{V}_{\theta \sharp}^{\perp}}\mu \right)\left\{ y \in B_E\left(x,2^{-k}\right) : d_E(\pi(x), \pi(y) )< \delta \right\}\right)^{p-1} \\ d\left( P_{\mathbb{V}_{\theta \sharp}^{\perp}}\mu_{F_k} \right)(x) \, d\theta 
\lesssim  \sum_{B \in \mathcal{B}} \int_0^{\pi} \int \\
 \left(\delta^{-1}\left(P_{\mathbb{V}_{\theta \sharp}^{\perp}}\mu \right)\left\{ y \in B_E\left(x,2^{-k}\right) : d_E(\pi(x), \pi(y) )< \delta \right\}\right)^{p-1}  d\left( P_{\mathbb{V}_{\theta \sharp}^{\perp}}\mu_B \right)(x) \, d\theta, \end{multline}
where $\mu_B$ is the restriction of $\mu_{F_k}$ to $B$. Let 
\begin{multline} \label{badballdefn} \mathcal{B}_b = \\
\left\{ B \in \mathcal{B} : \mathcal{H}^1\left\{ \theta \in [0, \pi) : P_{\mathbb{V}_{\theta \sharp}^{\perp}}\mu\left(P_{\mathbb{V}_{\theta}^{\perp}}(100B) \right) \geq c_{s,\mathbb{H}}(\mu) 2^{-k\left( s-1 - \epsilon \right) } \right\} \geq 2^{-k\eta} \right\}, \end{multline}
and let $\mathcal{B}_g = \mathcal{B} \setminus \mathcal{B}_b$.  Let $\mu_b = \sum_{B \in \mathcal{B}_b} \mu_B$, and $\mu_g = \sum_{B \in \mathcal{B}_g} \mu_B$. Then
\begin{multline} \label{goodbaddecomp} \eqref{sumoverB} \lesssim \int_0^{\pi} \int \left(\delta^{-1}\left(P_{\mathbb{V}_{\theta \sharp}^{\perp}}\mu \right)\left\{ y  : d_E(\pi(x), \pi(y) )< \delta \right\}\right)^{p-1}  d\left( P_{\mathbb{V}_{\theta \sharp}^{\perp}}\mu_b \right)(x) \, d\theta \\
+\int_0^{\pi} \int  \left(\delta^{-1}\left(P_{\mathbb{V}_{\theta \sharp}^{\perp}}\mu \right)\left\{ y \in B_E\left(x,2^{-k}\right) : d_E(\pi(x), \pi(y) )< \delta \right\}\right)^{p-1} \\
 d\left( P_{\mathbb{V}_{\theta \sharp}^{\perp}}\mu_g \right)(x) \, d\theta. \end{multline}
Suppose first that the term from $\mu_b$ dominates in \eqref{goodbaddecomp}. Then
\[ \eqref{sumoverB} \lesssim \int_0^{\pi} \int 
 \left(\delta^{-1}\left(\pi_{\sharp} P_{\mathbb{V}_{\theta \sharp}^{\perp}}\mu \right)\left\{ y  : d_E(x, y )< \delta \right\}\right)^{p-1}  d\left( \pi_{\sharp} P_{\mathbb{V}_{\theta \sharp}^{\perp}}\mu_b \right)(x) \, d\theta. \]
By Hölder's inequality, 
\begin{multline*} \eqref{sumoverB} \lesssim \left(\int_0^{\pi} \int 
 \left(\delta^{-1}\left(\pi_{\sharp} P_{\mathbb{V}_{\theta \sharp}^{\perp}}\mu \right)\left\{ y  : d_E(x, y )< \delta \right\}\right)^{p}  dx \, d\theta\right)^{1/p'} \\
\times \left(\int_0^{\pi} \int \left\lvert  \pi_{\sharp} P_{\mathbb{V}_{\theta \sharp}^{\perp}}\mu_b(x) \right\rvert^{p} \, dx \, d\theta\right)^{1/p}. \end{multline*}
The term $\delta^{-1}\left(\pi_{\sharp} P_{\mathbb{V}_{\theta \sharp}^{\perp}}\mu \right)\left\{ y  : d_E(x, y )< \delta \right\}$ is bounded by $M\pi_{\sharp} P_{\mathbb{V}_{\theta \sharp}^{\perp}}\mu(x)$, where $M$ is the Hardy-Littlewood maximal operator in one dimension. By the boundedness of the Hardy-Littlewood maximal operator on $L^p(\mathbb{R})$ applied to the first factor\footnote{Young's convolution inequality $\|f \ast g\|_p \leq \|f\|_p \|g\|_1$ could be used with $g = \delta^{-1} \chi_{(-\delta,\delta)}$ in place of the Hardy-Littlewood maximal inequality to avoid a constant that tends to $\infty$ as $p \to 1^+$, but using the maximal inequality is slightly cleaner.}, followed by an application of Theorem~\ref{lpestimate} with $R \sim 1$ to both factors, 
\begin{equation} \label{pause7000} \eqref{sumoverB} \lesssim \left(\mu(\mathbb{H}) c_{s-1,E}(\mu)^{p-1}\right)^{1/p'} \left(  \mu_b(\mathbb{H}) c_{s-1,E}(\mu)^{p-1}) \right)^{1/p}. \end{equation}
By the dimension comparison principle (\cite[Theorem~2.7]{BDFMT}, or more precisely \cite[Proposition~3.4]{BRS} from the proof of dimension comparison), 
\begin{equation} \label{comparisonfrostman} c_{s-1,E}(\mu) \lesssim c_{s,\mathbb{H}}(\mu). \end{equation}
 Lemma~\ref{projectionconjecture} implies that for $k$ sufficiently large,
\begin{equation} \label{quantprojbound} \mu_b(\mathbb{H}) = \sum_{B \in \mathcal{B}_b} \mu(B) \leq 2^{-k\eta} \mu(\mathbb{H}), \end{equation}
for $\eta>0$ sufficiently small depending only on $s$ and $\epsilon$. Substituting \eqref{comparisonfrostman} and \eqref{quantprojbound} into \eqref{pause7000} yields 
\begin{equation} \label{referencedagain} \eqref{sumoverB} \lesssim \mu(\mathbb{H}) c_{s,\mathbb{H}}(\mu)^{p-1} 2^{-k\eta/p}. \end{equation}
If $p>1$ is chosen sufficiently close to 1 (after $\eta$), this is stronger than \eqref{reusedbelow11}, so this proves the required inequality \eqref{reusedbelow11} in case the term from $\mu_b$ dominates in \eqref{goodbaddecomp}. 

Now suppose that the $\mu_g$ term dominates in \eqref{goodbaddecomp}.  Using the definition of pushforward, and then Fubini, 
\begin{multline} \label{unravel} \eqref{sumoverB} \lesssim \int \int_0^{\pi} \\
\left(\delta^{-1}\left(P_{\mathbb{V}_{\theta \sharp}^{\perp}}\mu \right)\left\{ y \in B_E\left(P_{\mathbb{V}_{\theta}^{\perp}}(x),2^{-k}\right) : d_E\left(\pi\left(P_{\mathbb{V}_{\theta}^{\perp}}(x)\right), \pi(y) \right)< \delta \right\}\right)^{p-1} \\
 d\theta \, d\mu_{g}(x). \end{multline}
After passing to a subset, it may be assumed that the balls $B \in \mathcal{B}_g$ are disjoint. For each $x$ in the support of $\mu_{g}$, choose a unique $B \in \mathcal{B}_g$ such that $x \in B$, and define 
\begin{equation} \label{badxdefn} \Theta_{b,x} = \left\{ \theta \in [0, \pi) : P_{\mathbb{V}_{\theta \sharp}^{\perp}}\mu\left(P_{\mathbb{V}_{\theta}^{\perp}}(10B) \right) \geq c_{s,\mathbb{H}}(\mu) 2^{-k\left( s-1-\epsilon \right) } \right\}, \end{equation}
and 
\begin{equation} \label{goodxdefn} \Theta_{g,x} = [0, \pi) \setminus \Theta_{b,x}. \end{equation} Then by \eqref{unravel},
\begin{multline} \label{thetadecomp} \eqref{sumoverB} \lesssim  \int \int_{\Theta_{b,x}}\\
 \left(\delta^{-1}\left(P_{\mathbb{V}_{\theta \sharp}^{\perp}}\mu \right)\left\{ y \in B_E\left(P_{\mathbb{V}_{\theta}^{\perp}}(x),2^{-k}\right) : d_E\left(\pi\left(P_{\mathbb{V}_{\theta}^{\perp}}(x)\right), \pi(y) \right)< \delta \right\}\right)^{p-1} \\
 \, d\theta \, d\mu_{g}(x) 
+ \int \int_{\Theta_{g,x}}\\
 \left(\delta^{-1}\left(P_{\mathbb{V}_{\theta \sharp}^{\perp}}\mu \right)\left\{ y \in B_E\left(P_{\mathbb{V}_{\theta}^{\perp}}(x),2^{-k}\right) : d_E\left(\pi\left(P_{\mathbb{V}_{\theta}^{\perp}}(x)\right), \pi(y) \right)< \delta \right\}\right)^{p-1} \\
 \, d\theta \, d\mu_{g}(x). \end{multline}
Consider the sub-case where the integral over $\Theta_{b,x}$ dominates the right-hand side of \eqref{thetadecomp}. Let $q >1$ be an exponent to be chosen. By the definition of $\mathcal{B}_g$ and $\mathcal{B}_b$ (see \eqref{badballdefn}), $\mathcal{H}^1\left( \Theta_{b,x} \right) \leq 2^{-k\eta}$ for each $x$ in the support of $\mu_{g}$. Hence, by Hölder's inequality, 
\begin{multline*} \eqref{sumoverB} \lesssim \Bigg(\int \int_0^{\pi} \\
\bigg(\delta^{-1}\left(P_{\mathbb{V}_{\theta \sharp}^{\perp}}\mu \right)\left\{ y \in B_E\left(P_{\mathbb{V}_{\theta}^{\perp}}(x),2^{-k}\right) : d_E\left(\pi\left(P_{\mathbb{V}_{\theta}^{\perp}}(x)\right), \pi(y) \right)< \delta \right\}\bigg)^{q(p-1)} \\
 \, d\theta \, d\mu_{g}(x)\Bigg)^{1/q} \left( \mu(\mathbb{H}) 2^{-k\eta} \right)^{1/q'}. \end{multline*}
Using Fubini and the definition of pushforward again, this can be relaxed to
\begin{multline*} \eqref{sumoverB} \lesssim \Bigg( \int_0^{\pi} \int \left(\delta^{-1}\left(P_{\mathbb{V}_{\theta \sharp}^{\perp}}\mu \right)\left\{ y  : d_E(\pi(x), \pi(y) )< \delta \right\}\right)^{q(p-1)} \\
 \, d\left(P_{\mathbb{V}_{\theta \sharp}^{\perp}}\mu\right)(x) \,  \, d\theta \Bigg)^{1/q} \left( \mu(\mathbb{H}) 2^{-k\eta} \right)^{1/q'}. \end{multline*}
If $p>1$ is sufficiently close to 1, and $q>1$ is defined such that $\widetilde{p} := q(p-1) + 1 = 4/3$, or equivalently $q = \frac{1}{3(p-1)}$, then by the boundedness of the Hardy-Littlewood maximal operator on $L^{\widetilde{p}}(\mathbb{R})$, and by Theorem~\ref{lpestimate} with $R \sim 1$ and with $\widetilde{p}=4/3$ in place of $p$,  
\begin{equation} \label{intheabove} \eqref{sumoverB} \lesssim  \left( \mu(\mathbb{H}) c_{s-1,E}(\mu)^{q(p-1)} \right)^{1/q}\left( \mu(\mathbb{H}) 2^{-k\eta} \right)^{1/q'}.\end{equation}
Using the dimension comparison inequality \eqref{comparisonfrostman}, and since $q' \to 1$ as $p \to 1^+$, \eqref{intheabove} will be stronger than \eqref{reusedbelow11} if $p$ is sufficiently close to 1, so this proves the required inequality \eqref{reusedbelow11} in the sub-case where the term from $\Theta_{b,x}$ dominates the right-hand side of \eqref{thetadecomp}.

It remains to consider the sub-case where the term from $\Theta_{g,x}$ dominates the right-hand side of \eqref{thetadecomp}. In this case, 
\begin{multline*}  \eqref{sumoverB} \lesssim  \sum_{B \in \mathcal{B}_g} \int \int_{\Theta_{g,x}} \\
\left(\delta^{-1}\left(P_{\mathbb{V}_{\theta \sharp}^{\perp}}\mu \right)\left\{ y \in B_E\left(P_{\mathbb{V}_{\theta}^{\perp}}(x),2^{-k}\right) : d_E\left(\pi\left(P_{\mathbb{V}_{\theta}^{\perp}}(x)\right), \pi(y) \right)< \delta \right\}\right)^{p-1} \\
 \, d\theta \, d\mu_{B}(x). \end{multline*}
By abbreviating $\Theta_{g,x} = \Theta_{g,B}$ when $x \in B$, and using Fubini and the definition of pushforward, this can be simplified to
\begin{multline} \label{pause721}  \eqref{sumoverB} \lesssim  \sum_{B \in \mathcal{B}_g}  \int_{\Theta_{g,B}} \int \\
 \left(\delta^{-1}\left(P_{\mathbb{V}_{\theta \sharp}^{\perp}}\mu \right)\left\{ y \in B_E\left(x,2^{-k}\right) : d_E(\pi(x), \pi(y) )< \delta \right\}\right)^{p-1}  \, d\left(P_{\mathbb{V}_{\theta \sharp}^{\perp}}\mu_{B}\right)(x) \, d\theta. \end{multline}
For each $B \in \mathcal{B}_g$ and each $\theta \in \Theta_{g,B}$, the $\mu$ in the integrand can be replaced by 
\[ \mu \chi_{T_{B, \theta}} =: \mu_{T_{B, \theta}}, \qquad  T_{B, \theta}  = P_{\mathbb{V}_{\theta}^{\perp}}^{-1}\left(P_{\mathbb{V}_{\theta}^{\perp}}(10B)\right). \]
An important inequality will be that for any $B \in \mathcal{B}_g$ and any $\theta \in \Theta_{g,B}$, 
\begin{equation} \label{bushbound1} \mu_{T_{B, \theta} }(\mathbb{H}) = \mu\left(T_{B, \theta}\right) \leq  c_{s,\mathbb{H}}(\mu) 2^{-k\left( s-1 -\epsilon \right)}, \end{equation}
which follows from the definition of $T_{B, \theta}$ and $\Theta_{g,B}$ when $B \in \mathcal{B}_g$ (see \eqref{badxdefn} and \eqref{goodxdefn}). Therefore, \eqref{pause721} can be relaxed to 
\begin{multline*}   \eqref{sumoverB} \lesssim  \sum_{B \in \mathcal{B}_g}  \int_{\Theta_{g,B}} \int \left(\delta^{-1}\left(P_{\mathbb{V}_{\theta \sharp}^{\perp}}\mu_{T_{B, \theta}} \right)\left\{ y : d_E(\pi(x), \pi(y) )< \delta \right\}\right)^{p-1} \\
 \, d\left(P_{\mathbb{V}_{\theta \sharp}^{\perp}}\mu_{B}\right)(x) \, d\theta . \end{multline*}
For each $B$, decompose $\mu_B$ into ``high'' and ``low'' frequencies:
\[ \mu_B =  \mu_B \ast \widecheck{\psi_l} +  \mu_B \ast \widecheck{\psi_h}, \]
where $\psi_l$ is a smooth bump function on $|\xi| \lesssim 2^{k}$,  $\psi_h$ is a smooth bump on $|\xi| \gtrsim 2^k$. Then 
\begin{multline} \label{freqdecomp} \eqref{sumoverB} \lesssim  \sum_{B \in \mathcal{B}_g} \bigg\lvert \int_{\Theta_{g,B}} \int \left(\delta^{-1}\left(P_{\mathbb{V}_{\theta \sharp}^{\perp}}\mu_{T_{B, \theta} } \right)\left\{ y : d_E(\pi(x), \pi(y) )< \delta \right\}\right)^{p-1} \\
 d\left( P_{\mathbb{V}_{\theta \sharp}^{\perp}}\left(\mu_B \ast\widecheck{\psi_l}\right)  \right)(x) \, d\theta \bigg\rvert \\
+ \sum_{B \in \mathcal{B}_g} \bigg\lvert \int_{\Theta_{g,B}} \int 
 \left(\delta^{-1}\left(P_{\mathbb{V}_{\theta \sharp}^{\perp}}\mu_{T_{B, \theta} } \right)\left\{ y : d_E(\pi(x), \pi(y) )< \delta \right\}\right)^{p-1} \\
 d\left( P_{\mathbb{V}_{\theta \sharp}^{\perp}}\left(\mu_B \ast\widecheck{\psi_h}\right) \right)(x) \, d\theta \bigg\rvert. \end{multline} 
Suppose that the first term in \eqref{freqdecomp} dominates (low frequency case), and therefore 
\begin{multline} \label{smoothedversion} \eqref{sumoverB} \lesssim  \sum_{B \in \mathcal{B}_g} \int_{\Theta_{g,B}} \int \left(\delta^{-1}\left(P_{\mathbb{V}_{\theta \sharp}^{\perp}}\mu_{T_{B, \theta} } \right)\left\{ y  : d_E(\pi(x), \pi(y) )< \delta \right\}\right)^{p-1} \\
\, d\left( P_{\mathbb{V}_{\theta \sharp}^{\perp}}\mu_{B,k}  \right)(x) \, d\theta, \end{multline}
where $\mu_{B,k} = \mu_B \ast \phi_k$, with $\phi_k$ a non-negative smooth bump function satisfying $\phi_k \lesssim 2^{3k}$, with $\phi_k$ rapidly decaying outside the Euclidean ball $B_{E}\left(0, 2^{-k}\right)$. Each measure $\mu_{B,k}$ is rapidly decaying outside $B$. By the definition of pushforward applied to \eqref{smoothedversion}, 
\begin{multline*} \eqref{sumoverB} \lesssim  \sum_{B \in \mathcal{B}_g} \int_{\Theta_{g,B}} \int \left(\delta^{-1}\left(\pi_{\sharp} P_{\mathbb{V}_{\theta \sharp}^{\perp}}\mu_{T_{B, \theta} } \right)\left\{ y  : d_E(x, y)< \delta \right\}\right)^{p-1} \\
\, d\left(\pi_{\sharp} P_{\mathbb{V}_{\theta \sharp}^{\perp}}\mu_{B,k}  \right)(x) \, d\theta. \end{multline*}
By Hölder's inequality with $q= \frac{1}{p-1}$, folowed by Young's convolution inequality $\lVert f \ast g \rVert_1 \leq \lVert f \rVert_1 \lVert g \rVert_1$ with $g= \delta^{-1} \chi_{(-\delta,\delta)}$ (or just Fubini), this becomes 
\begin{equation} \label{referencedagain2} \eqref{sumoverB} \lesssim 2^{k(p-1)O(\epsilon)}\sum_{B \in \mathcal{B}_g} \left( \int_{\Theta_{g,B}} \mu(T_{B, \theta}) \, d\theta\right)^{1/q} \left( \mu(B)^{q'} 2^{k(q'-1)} \right)^{1/q'}, \end{equation}
where, for the second factor, the trivial $L^{q'}$ inequality for the projection was used since $\mu_{B,k}$ can be treated as constant (more precisely $\mu_{B,k} \lesssim 2^{3k} 2^{kO(\epsilon) } \mu(B) \chi_{B}$, where $\chi_B$ is $\lesssim 1$ on $B$ and rapidly decaying outside $B$). Applying \eqref{bushbound1} to the above gives 
\[ \eqref{sumoverB} \lesssim  2^{k(p-1)O(\epsilon)} \sum_{B \in \mathcal{B}_g} \left( c_{s,\mathbb{H}}(\mu) 2^{-k \left( s-1 - \epsilon \right)} \right)^{1/q}  \left( \mu(B)^{q'} 2^{k(q'-1)} \right)^{1/q'}. \] 
This simplifies to 
\[ \eqref{sumoverB} \lesssim   \mu(\mathbb{H})  c_{s,\mathbb{H}}(\mu)^{p-1} 2^{-k\left( s-2 - O(\epsilon) \right)(p-1)  }. \]
This verifies the required inequality \eqref{reusedbelow11} in the case where the first term in \eqref{freqdecomp} dominates. 

 If the second term dominates in \eqref{freqdecomp} (high frequency case), then by the definition of pushforward, followed by Hölder's inequality with $q = \frac{1}{p-1}$, 
\begin{multline} \label{reusedagain} \eqref{sumoverB} \lesssim  \sum_{B \in \mathcal{B}_g} \left( \int_{\Theta_{g,B}} \int  \delta^{-1}\left(\pi_{\sharp} P_{\mathbb{V}_{\theta \sharp}^{\perp}}\mu_{T_{B, \theta} } \right)\left\{ y  : d_E(x,y )< \delta \right\} \, dx \, d\theta\right)^{1/q} \\
\times  \left( \int_0^{\pi} \int \left\lvert \pi_{\sharp} P_{\mathbb{V}_{\theta \sharp}^{\perp}}(\mu_B \ast\widecheck{\psi_h})  \right\rvert^{q'} \, \mathcal{H}^1_E \, d\theta\right)^{1/{q'}}. \end{multline}
By Young's convolution inequality $\lVert f \ast g \rVert_1 \leq \lVert f \rVert_1 \lVert g \rVert_1$ (or just Fubini) applied to the first factor above for each $j$, and the second part of Theorem~\ref{lpestimate} with $R\sim 2^k$ applied to the second factor (using \eqref{upperdensity2}), this gives, for small $\epsilon >0$,
\begin{equation} \label{reusedbelow21} \eqref{sumoverB} \lesssim  \\
     \sum_{B \in \mathcal{B}_g} \left( \int_{\Theta_{g,B}} \mu(T_{B, \theta}) \, d\theta \right)^{1/q} \left(  \mu(B) 2^{-k(s_E-1-O(\epsilon))(q'-1)}\right)^{1/q'}. \end{equation}
The dimension comparison inequality \eqref{comparisonfrostman}, together with \eqref{lowerdensity2} and \eqref{upperdensity2}, was used to obtain that $s_E \geq s -1-O(\epsilon) >1$. Therefore, applying \eqref{bushbound1} to the first factor above gives
\begin{equation} \label{notfinished} \eqref{sumoverB} \lesssim  c_{s,\mathbb{H}}(\mu)^{p-1} 2^{-k(s-1-\epsilon)/q} 2^{-k(s_E-1-O(\epsilon))/q} \sum_{B \in \mathcal{B}_g} \mu(B)^{1/q'}   . \end{equation}
By the lower density inequality \eqref{lowerdensity2}, and the property that each ball $B \in \mathcal{B}$ intersects the support of $\mu_{F_k}$, each $B \in \mathcal{B}$ satisfies $\mu(B) \gtrsim 2^{-k(s_E+O(\epsilon))}$. Hence
\[ |\mathcal{B}_g| \leq |\mathcal{B} | \lesssim \mu(\mathbb{H}) 2^{k(s_E+\epsilon ) }. \]
Therefore, applying Hölder's inequality to the sum in \eqref{notfinished} yields 
\[ \eqref{sumoverB} \lesssim c_{s,\mathbb{H}}(\mu)^{p-1} \mu(\mathbb{H})  2^{-k(s-2-O(\epsilon))(p-1)}. \]
This verifies \eqref{reusedbelow11}, so this proves the required inequality \eqref{reusedbelow11} in the case where the sum over $j$ dominates in \eqref{freqdecomp}. This finishes the proof of the second statement \eqref{secondclaim} in the theorem.

To prove the first statement \eqref{firstclaim}, the only adjustments needed are that replacing $\mu$ by $\mu_F$ is no longer necessary since the Euclidean lower density assumption in \eqref{lowerdensity2} is automatically satisfied, and the decomposition of $F$ into $F_k$ is no longer necessary. Since the lower density assumption is satisfied, the outer $\liminf$ in the integral in \eqref{desiredinequality} can be replaced by a $\limsup$, which can be similarly bounded by a sum over a tail of terms exceeding $K$ as in \eqref{sumbound}, for arbitrarily large $K$. Bounding the resulting series is virtually identical to the rest of the proof of \eqref{secondclaim}.   \end{proof}

The remainder of the proof of Theorem~\ref{etccorollary1} is given below. 

\begin{proof}[Proof of Theorem~\ref{etccorollary1}]  Let $A \subseteq \mathbb{H}$ be a Borel set with $2 < \dim A < 3$. Let $t = \dim A$ and let $\epsilon >0$ with $\epsilon \ll t-2$. By the second part of Theorem~\ref{intersectiontheorem2}, $\mathcal{H}^1_E(F_{\theta}) >0$ for a.e.~$\theta \in [0, \pi)$, where 
\[ F_{\theta}= \left\{ \lambda   \in \mathbb{R} : \dim_P\left( P_{\mathbb{V}_{\theta}^{\perp}}(A)  \cap \pi^{-1}(\lambda) \right) > t-2 - \epsilon \right\}. \]
 Let $\theta \in [0, \pi)$ be such that $\mathcal{H}^1_E(F_{\theta})>0$. Let $\{E_{\theta,i}\}_i$ be a countable covering of $P_{\mathbb{V}_{\theta}^{\perp}}(A)$ by compact subsets of $\mathbb{V}_{\theta}^{\perp}$. By countable stability of the packing dimension, there exists $i$ such that $\mathcal{H}^1_E(F_{\theta,i}) >0$, where 
\[ F_{\theta,i}= \left\{ \lambda \in \mathbb{R} : \dim_P\left( E_{i,\theta}  \cap \pi^{-1}(\lambda) \right) > t-2 - \epsilon \right\}. \]
Let $\delta>0$ be small. For each $\lambda \in F_{\theta,i}$, since the upper Minkowski dimension is greater than or equal to the packing dimension, there exists a collection $\{B(x_{\theta,\lambda,j}, r_{\theta,\lambda}) \}_{j=1}^{N(\theta, \lambda)}$ of disjoint intervals of equal dyadic radii $r_{\theta, \lambda} < \delta$, centred at points $x_{\theta,\lambda,j} \in E_{\theta,i}  \cap \pi^{-1}(\lambda)$, such that 
\[ \sum_{j=1}^{N(\theta, \lambda)} r_{\theta, \lambda}^{t-2-\epsilon} >1. \] 
By dyadic pigeonholing, there is a fixed integer $k = k(\theta)$ independent of $\lambda$ and a set $F_{\theta,i,k} \subseteq F_{\theta,i}$ with $\mathcal{H}^1_E(F_{\theta,i,k}) \gtrsim k^{-2} \mathcal{H}^1_E(F_{\theta,i})$ such that $r_{\theta,\lambda} = 2^{-k}$ for all $\lambda \in F_{\theta,i,k}$. Let $\Lambda$ be a maximal $2^{-2k}$-separated subset of $F_{\theta,i,k}$. Then the cardinality of $\Lambda$ satisfies $|\Lambda| \gtrsim 2^{2k} k^{-2} \mathcal{H}^1_E(F_{\theta,i})$. Therefore, the family $\mathcal{R}(\theta) = \{ B_{\mathbb{H}}(x_{\theta,\lambda,j}, 2^{-k} ) : \lambda \in \Lambda, \quad 1 \leq j \leq  N(\theta, \lambda) \}$ is a disjoint family of Korányi balls (which can be thought of as $\sim 2^{-k} \times 2^{-2k}$ rectangles when intersected with $\mathbb{V}_{\theta}^{\perp}$), centred at points in $E_{\theta,i}$, of radii $2^{-k} < \delta$, such that,
\[  \left\lvert \mathcal{R}_{\theta} \right\rvert \geq 2^{k(t-2\epsilon)} \left(2^{k\epsilon} k^{-2}\mathcal{H}^1(F_{\theta,i}) \right). \] 
This proves that for any sufficiently small $\delta>0$, there is some $2^{-k} < \delta$ such that the number of Korányi balls of radius $2^{-k}$ required to cover $E_{\theta,i}$ is at least $2^{k(t-2\epsilon)}$. By definition, this shows that the upper Minkowski dimension of $E_{\theta,i}$ exceeds $t-2\epsilon$. Since $\{E_{\theta,i}\}_i$ is an arbitrary countable covering of $P_{\mathbb{V}_{\theta}^{\perp}}(A)$ by compact subsets of $\mathbb{V}_{\theta}^{\perp}$, this proves that $\dim_P P_{\mathbb{V}_{\theta}^{\perp}}(A) \geq t-2\epsilon$ for a.e.~$\theta \in [0, \pi)$. Since this holds for arbitrarily small $\epsilon >0$, it follows that $\dim_P P_{\mathbb{V}_{\theta}^{\perp}}(A) \geq t$ for a.e.~$\theta \in [0, \pi)$. \end{proof}

The rest of the proof of Proposition~\ref{etccorollary} is given below. 

\begin{proof}[Proof of Proposition~\ref{etccorollary}] Let $A \subseteq \mathbb{H}$ be a Borel set. For any $\theta \in [0, \pi)$ and $2 < t \leq 3$, it will be shown that
\begin{equation} \label{claimedinequality}  \mathcal{H}^{t}_{\mathbb{H}}\left(P_{\mathbb{V}_{\theta}^{\perp}}(A)\right) \gtrsim \int_{\mathbb{R}}^{*} \mathcal{H}^{t-2}\left( \pi^{-1}(\lambda) \cap P_{\mathbb{V}_{\theta}^{\perp}}(A) \right) \, d\lambda, \end{equation}
where $\int^* f$ refers to the upper integral of $f$, defined as the infimum of $\int g$ over Lebesgue measurable functions $g \geq f$. The inequality \eqref{claimedinequality} follows from the same argument as in \cite[Theorem~7.7]{mattila}, especially considering the projections $\pi$ are Lipschitz with respect to the Korányi metric when restricted to vertical planes, but the details are included below. By definition, for any $\theta \in [0, \pi)$,
\begin{equation} \label{pause230} \int_{\mathbb{R}}^{*} \mathcal{H}^{t-2}\left( \pi^{-1}(\lambda) \cap P_{\mathbb{V}_{\theta}^{\perp}}(A) \right) \, d\lambda =\int_{\mathbb{R}}^{*}  \liminf_{k \to \infty}  \mathcal{H}^{t-2}_{2^{-k}}\left( \pi^{-1}(\lambda) \cap P_{\mathbb{V}_{\theta}^{\perp}}(A) \right) \, d\lambda. \end{equation}
Fix a large integer $k$, and let $\{B_{\mathbb{H}} (z_{j,k}, r_{j,k})\}_j$ be a covering of $P_{\mathbb{V}_{\theta}^{\perp}}(A)$ by Korányi balls of radius at most $2^{-k}$ and centres in $\mathbb{V}_{\theta}^{\perp}$, such that 
\[ \sum_{j} r_{j,k}^t \leq \mathcal{H}^t_{2^{-k}}\left(P_{\mathbb{V}_{\theta}^{\perp}}(A)\right) + \frac{1}{k}. \]
Then
\begin{multline} \label{subadditivity} \int_{\mathbb{R}}^{*}  \liminf_{k \to \infty}  \mathcal{H}^{t-2}_{2^{-k}}\left( \pi^{-1}(\lambda) \cap P_{\mathbb{V}_{\theta}^{\perp}}(A) \right) \, d\lambda \\
\lesssim \int_{\mathbb{R}}  \liminf_{k \to \infty} \left( \sum_j \diam\left( \pi^{-1}(\lambda) \cap \mathbb{V}_{\theta}^{\perp} \cap B_{\mathbb{H}} (z_{j,k}, r_{j,k} ) \right)^{t-2} \right) \, d\lambda. \end{multline}
For each $j$ and $k$, let 
\[ F_{j,k} = \left\{ \lambda \in \mathbb{R} : \pi^{-1}(\lambda) \cap  \mathbb{V}_{\theta}^{\perp} \cap  B_{\mathbb{H}}(z_{j,k}, r_{j,k}) \neq \emptyset \right\}. \]
Since the upper integral has been replaced by a standard integral in \eqref{subadditivity}, Fatou's lemma and the monotone convergence theorem can be used to obtain
\[ \eqref{pause230} \leq  \liminf_{k \to \infty} \sum_j \int_{F_{j,k}} r_{j,k}^{t-2} \, d\lambda. \]
Each Korányi ball $B_{j,k}(z_{j,k}, r_{j,k})$ intersected with $\mathbb{V}_{\theta}^{\perp}$ is contained in a rectangle of dimensions $2r_{j,k}  \times \frac{1}{2}r_{j,k}^2$, with the last coordinate in the vertical direction, and therefore $F_{j,k}$ is contained in an interval of length $\frac{1}{2}r_{j,k}^2$. Hence,
\[ \eqref{pause230} \lesssim \liminf_{k \to \infty} \sum_j r_{j,k}^{t-2} r_{j,k}^{2} \leq \lim_{k \to \infty}\left(\mathcal{H}^t_{2^{-k}}(P_{\mathbb{V}_{\theta}^{\perp}}(A) ) + \frac{1}{k}\right) = \mathcal{H}^t_{\mathbb{H}}(P_{\mathbb{V}_{\theta}^{\perp}}(A) ). \] 
This verifies \eqref{claimedinequality}.

Now let $\mu$ be a Borel meausure on $\mathbb{H}$ which is Euclidean Ahlfors-regular, let $t = \dim^*_{\mathbb{H}} \mu$, and assume that $2 < t < 3$ (the other cases are simpler and follow from known results, and can also be deduced from Theorem~\ref{lpestimate}). Then, for any $\epsilon >0$ (with $\epsilon \ll t-2$), there is a Borel set $A \subseteq \supp \mu$ with $\mu(A)>0$ such that $\mu(B_{\mathbb{H}}(x,r)) \lesssim r^{t-\epsilon}$ for all $x \in A$ and $r>0$.  Let $F \subseteq [0, \pi)$ be the set of $\theta \in [0, \pi)$ for which $\mathcal{H}^{t-2\epsilon}_{\mathbb{H}}\left(P_{\mathbb{V}_{\theta}^{\perp}}(A)\right) = 0$.  By \eqref{claimedinequality} with $t-2\epsilon$ in place of $t$,
\begin{equation} \label{pause25} 0 = \int_F \mathcal{H}^{t-2\epsilon}_{\mathbb{H}}\left(P_{\mathbb{V}_{\theta}^{\perp}}(A)\right) \, d \theta \geq \int_F \int_{\mathbb{R}}^{*} \mathcal{H}^{t-2-2\epsilon}\left( \pi^{-1}(\lambda) \cap P_{\mathbb{V}_{\theta}^{\perp}}(A) \right) \, d\lambda \, d\theta. \end{equation}

Since $t-2-2\epsilon <t-2-\epsilon$, the first part of Theorem~\ref{intersectiontheorem2} implies that for a.e.~$\theta \in [0, \pi)$, there is a positive length set of $\lambda$ such that 
\[ \mathcal{H}^{t-2-2\epsilon}_{\mathbb{H}}\left( \pi^{-1}(\lambda) \cap P_{\mathbb{V}_{\theta}^{\perp}}(A) \right) = \infty, \]
 and thus for a.e.~$\theta \in [0, \pi)$, 
\[ \int_{\mathbb{R}}^{*} \mathcal{H}^{t-2-2\epsilon}\left( \pi^{-1}(\lambda) \cap P_{\mathbb{V}_{\theta}^{\perp}}(A) \right) \, d\lambda  =  \infty. \]
Substituting into \eqref{pause25} yields that $\mathcal{H}^1(F) = 0$, or equivalently $\mathcal{H}^{t-2\epsilon}_{\mathbb{H}}\left(P_{\mathbb{V}_{\theta}^{\perp}}(A)\right) >0$ for a.e.~$\theta \in [0, \pi)$. Therefore, $\mathcal{H}^{t-2\epsilon}_{\mathbb{H}}\left(P_{\mathbb{V}_{\theta}^{\perp}}(\supp \mu)\right) >0$ for a.e.~$\theta \in [0, \pi)$. Since, in the latter statement, $\epsilon>0$ can be taken arbitrarily small, it follows that $\dim_{\mathbb{H}}\left(P_{\mathbb{V}_{\theta}^{\perp}}(\supp \mu)\right) \geq t$ for a.e.~$\theta \in [0, \pi)$. \end{proof} 

\section{Hausdorff dimension of projections}

The following lemma roughly says that if the projections of a measure satisfy an $s$-Frostman condition on average, then the projected supports of the measure almost surely have dimension at least $s$.
\begin{lemma} \label{averaginglemma} Let $\sum_{k=1}^{\infty} a_k$ be any convergent series of non-negative real numbers. Let $\Theta$ be a compact topological space and let $X$ be a separable metric space. Suppose that $(\theta, x) \mapsto \pi_{\theta}(x)$ is a continuous function from $\Theta \times X$ into $X$.  Let $s\geq 0$, and let $\nu$ be a Borel measure on $\Theta$. Suppose that $\mu$ is a compactly supported nonzero finite Borel measure on $X$, and suppose that there exists $K\geq 1$ such that for all $k \geq K$, 
\begin{equation} \label{assumedbound} \int \left( \pi_{\theta \sharp} \mu \right)\left( \bigcup_{D \in \mathbb{D}_{\theta}} D \right) \, d\nu(\theta) \leq a_k, \end{equation}
whenever each $\mathbb{D}_{\theta}$ is a set of radius $2^{-k}$ balls in $X$ of cardinality $\left\lvert \mathbb{D}_{\theta} \right\rvert \leq 2^{ks}$ such that the integrand of \eqref{assumedbound} is Borel measurable in $\theta$. Then $\dim\left( \pi_{\theta} (\supp\mu)\right) \geq s$ for $\nu$-a.e.~$\theta \in \Theta$. 
\end{lemma}
\begin{remark} The proof shows that $\dim\left( \pi_{\theta} (\supp\mu)\right) \geq s$ could be upgraded to $\mathcal{H}^s(\pi_{\theta} (\supp\mu)) = \infty$ where $\mathcal{H}^s$ refers to the Hausdorff measure on $X$, but this more precise version will not be needed. \end{remark}
\begin{proof}[Proof of Lemma~\ref{averaginglemma}] The separability assumption on $X$ implies that $\mu(X \setminus \supp \mu) = 0$, which is used below.

Let $E$ be the set of $\theta \in \Theta$ such that $\dim\left( \pi_{\theta} (\supp\mu)\right) < s$. It suffices to show that $\int^* \chi_E \, d\nu = 0$, where the upper integral $\int^* f$ of $f$ is the infimum over the integrals $\int g$ of measurable functions $g \geq f$. Let $N \geq K$. For each $\theta \in \Theta$, let $\mathbb{D}_{\theta}$ be a collection of open balls of dyadic radii at most $2^{-N}$,
which covers $\pi_{\theta}(\supp\mu)=\supp \pi_{\theta \sharp}\mu$ when $\theta \in E$, such that
\begin{equation} \label{hdimensiondefn}  \sum_{D \in \mathbb{D}_{\theta}} \diam(D)^{s} \leq 1, \end{equation} 
and such that the integrands of the right-hand sides of \eqref{measurable} and \eqref{nubound} below are Borel measurable in $\theta$ for any $k$. This can be arranged by using continuity and compactness to make the sets $\mathbb{D}_{\theta}$ piecewise constant over a finite partition of $\Theta$ into Borel sets. Then 
\begin{equation} \label{measurable} \int^* \chi_E \, d\nu \leq \frac{1}{\mu(X)} \int \left(\pi_{\theta \sharp}\mu\right)\left( \bigcup_{D \in \mathbb{D}_{\theta}} D \right) \, d\nu(\theta). \end{equation}
Let $\mathbb{D}_{\theta,k} = \{ D \in \mathbb{D}_{\theta} : \diam(D)/2 = 2^k \}$. Then \eqref{measurable} becomes
\begin{equation} \label{nubound} \int^* \chi_E \, d\nu \leq \frac{1}{\mu(X)} \sum_{k =N}^{\infty} \int \left(\pi_{\theta \sharp}\mu\right)\left( \bigcup_{D \in \mathbb{D}_{\theta,k}}D \right) \, d\nu(\theta). \end{equation}
By \eqref{hdimensiondefn}, $\left\lvert \mathbb{D}_{\theta,k} \right\rvert \leq 2^{ks}$ for all $\theta \in E$. Therefore, applying the assumed \eqref{assumedbound} to \eqref{nubound} gives 
\[ \int^* \chi_E \, d\nu \leq  \frac{1}{\mu(X)} \sum_{k=N}^{\infty} a_k. \]
Letting $N \to \infty$ gives $\int^* \chi_E \, d\nu=0$.    \end{proof}

For a set $X$ in Euclidean space and a dyadic number $\Delta$, let $\mathcal{D}_{\Delta}(X)$ denote the set of semi-open dyadic cubes intersecting $X$ (open on left, closed on right). The dyadic notation and the concept of uniformity (used below) are taken from \cite{orponenshmerkinabc}. 

\begin{proof}[Proof of Theorem~\ref{etccorollary2}] Let $\mu$ be a Borel measure supported in $\{(z,t) \in \mathbb{H} : |z| \sim 1\}$, with support contained in a semi-open dyadic cube of side length 1, such that $c_{t, \mathbb{H}}(\mu) \leq 1$ where $2 < t < 3$. Let 
\[ g(t) = \max\left\{ \frac{(t-1)(t-2)}{t^2-4t+6}, t -2-  \frac{(t-1)(t-2)(3-t)}{-3t^2+14t-14}, 2t-5 \right\}. \]
By Frostman's lemma and Lemma~\ref{averaginglemma}, it suffices to show that, for any positive $\varepsilon$ sufficiently small, for any $k$ sufficiently large depending on $\varepsilon$, if for every $\theta \in [0, \pi)$, $\{ \mathbb{D}_{\theta}\}$ is a family of $2^{-k} \times 2^{-2k}$ axis parallel rectangles in $\mathbb{V}_{\theta}^{\perp}$ of cardinality $\left\lvert \mathbb{D}_{\theta} \right\rvert \leq 2^{k\left(2+g(t) - \sqrt{\varepsilon} \right)}$, then 
\begin{equation} \label{goal} \int_0^{\pi} \int_{\bigcup_{D \in \mathbb{D}_{\theta}} D } \, d\left( \mathcal{P}_{\mathbb{V}_{\theta \sharp}^{\perp}}\mu\right) \, d\theta \leq 2^{-k \varepsilon}, \end{equation}
whenever the inner integral is Borel measurable in $\theta$. To avoid measurability technicalities below, it may be further assumed after identifying $\mathbb{V}_{\theta}^{\perp}$ with $\mathbb{R}^2$ that $\mathbb{D}_{\theta}$ is piecewise constant in $\theta$ over a partition of $[0, \pi]$ into intervals of length $2^{-2k}$. Let $\Delta$ be a dyadic number with $\Delta \sim 2^{-k\varepsilon}$, where $\varepsilon$ will implicitly be taken sufficiently small  in the argument below. Let $L$ be the largest positive integer such that $\Delta^L > 2^{-k}$. After pigeonholing \eqref{goal} (losing a factor of $\log(2^k)^{O(1)}$), it may be assumed that for any $1 \leq m \leq L$, the quantity $\mu(B)$ is dyadically constant\footnote{The need for this pigeonholing step seems to be an obstruction to an analogous intersection theorem for the Hausdorff dimension.} over $B \in \mathcal{D}_{\Delta^m}(\supp \mu)$, and over $B \in \mathcal{D}_{2^{-k}}(\supp \mu)$. It is straightforward to bound the left-hand side of \eqref{goal} over the rectangles $D$ with $\left( P_{\mathbb{V}_{\theta \sharp}^{\perp}}\mu\right)(D) \leq  2^{-k\left(2+g(t)\right)}2^{k\sqrt{\varepsilon}/2}$, so it may be assumed that 
\[ \left( P_{\mathbb{V}_{\theta \sharp}^{\perp}}\mu\right)(D) >  2^{-k\left(2+g(t)\right)}2^{k\sqrt{\varepsilon}/2}, \quad  \forall \, D \in \mathbb{D}_{\theta} \quad \forall \, \theta \in [0, \pi). \]
 Therefore, by taking $\delta=2^{-2k}$, the left-hand side of \eqref{goal} satisfies 
\begin{multline*} \eqref{goal} \lesssim 2^{k \left(g(t) - \sqrt{\varepsilon}/2\right)(p-1)} \int_0^{\pi} \int 
 \\\left(\delta^{-1} \left(P_{\mathbb{V}_{\theta \sharp}^{\perp}}\mu \right)\left\{ y \in B_E\left(x,2^{-k}\right) : d_E(\pi(x), \pi(y) )< \delta \right\}\right)^{p-1}  
 d\left( P_{\mathbb{V}_{\theta \sharp}^{\perp}}\mu \right)(x) \, d\theta, \end{multline*}
where $p>1$ will be chosen sufficiently close to 1 below. By the above, to prove \eqref{goal} it suffices to show that for $\varepsilon$ sufficiently small,
 \begin{multline} \label{reusedbelow12} \int_0^{\pi} \int 
 \left(\delta^{-1} \left(P_{\mathbb{V}_{\theta \sharp}^{\perp}}\mu \right)\left\{ y \in B_E\left(x,2^{-k}\right) : d_E(\pi(x), \pi(y) )< \delta \right\}\right)^{p-1}  \\
 d\left( P_{\mathbb{V}_{\theta \sharp}^{\perp}}\mu \right)(x) \, d\theta \lesssim  2^{-k(p-1)\left(g(t) - O(\varepsilon)\right)}. \end{multline}
Denote $s_E = t-1$. Let $\gamma_L= \min\{4-t,x\}$ where $x$ solves the equation 
\begin{equation} \label{gammaLdefn2} \frac{t-2}{2}\left( 1+ \frac{x}{t-1}\right)  = t-2 - x \min\left\{  \frac{1}{4-t} + \frac{1}{t-1}-1, \frac{3-t}{4-t} \right\}. \end{equation}
A formula for $\gamma_L$ is 
\begin{multline} \label{gammaLdefn} \gamma_{L} = \min\left\{\max\left\{\frac{(t-1)(t-2)(4-t)}{t^2-4t+6}, \frac{(t-1)(t-2)(4-t)}{-3t^2+14t - 14}\right\},4-t \right\} \\
= \begin{cases} \frac{(t-1)(t-2)(4-t)}{t^2-4t+6} & 2 < t < 5/2 \\
\frac{(t-1)(t-2)(4-t)}{-3t^2+14t - 14} & 5/2 \leq t < \frac{1}{8}\left( 17 + \sqrt{33}\right) \\
4-t & \frac{1}{8}\left( 17 + \sqrt{33}\right) \leq t < 3. \end{cases} \end{multline}
Let $\Gamma_L = L \gamma_L$, and define $\Gamma_{L-1}$ by 
\[ (L-1)(s_E-1) + \Gamma_{L-1} = L(s_E-1) +\Gamma_L\left( 1- \frac{1}{s_E} \right), \]
or equivalently
\[ \Gamma_{L-1} = s_E-1   + \Gamma_L\left( 1- \frac{1}{s_E} \right). \]
Define $\gamma_{L-1}$ by $\Gamma_{L-1} = (L-1) \gamma_{L-1}$, so that
\begin{equation} \label{gammaL1defn} \gamma_{L-1} = \frac{1}{L-1} \left(s_E-1   + L\gamma_L\left( 1- \frac{1}{s_E} \right)\right). \end{equation}
 Define 
$\gamma_1 \geq \gamma_2 \geq \dotsb \geq \gamma_{L-1}$ recursively by
\begin{equation} \label{recursion} \Gamma_m = \Gamma_{m+1}+ s_E-1, \qquad \Gamma_m = m \gamma_m, \qquad 1 \leq m \leq L-2. \end{equation}
Equivalently, 
\[ \Gamma_m = (L-m)(s_E-1) + \frac{\Gamma_L(s_E-1)}{s_E}, \qquad 1 \leq m \leq L-1, \]
or
\begin{equation} \label{gammadefn} \gamma_m = \frac{ (L-m)(s_E-1)}{m} + \frac{L \gamma_L(s_E-1)}{ms_E}, \qquad 1 \leq m \leq L-1. \end{equation}

Denote $\Delta_m = \Delta^m$. Let $l$ be the largest integer $l \in [1,L]$ with $\left\lvert \mathcal{D}_{\Delta_l}(\supp \mu) \right\rvert \leq \Delta_l^{-(s_E+\gamma_l)}$. This is well-defined provided $\varepsilon \ll s_E-1$, since then $L \gtrsim 1/\varepsilon$, which yields $\gamma_1 \gg 3$, so $s_E + \gamma_{1} \geq 3$, and the support of $\mu$ being contained in a semi-open dyadic cube in $\mathbb{R}^3$ of side length 1 then implies that $\left\lvert \mathcal{D}_{\Delta_1}(\supp \mu) \right\rvert \leq \Delta_1^{-(s_E+\gamma_1)}$. Therefore $l$ is well-defined. 
Moreover, $l \neq L-1$ if $\varepsilon$ is small enough, since $\gamma_L - \gamma_{L-1} \geq (1/2) (\gamma_L/s_E)$ if $\varepsilon$ is small enough. 

By similar reasoning to the above, and by maximality of $l$,  $s_E +\gamma_{l+1} < 3$ if $l < L$. 
This implies that if $l < L-1$
\begin{equation} \label{lbound} l +1 > 
L\left(\frac{s_E-1}{2}\right)\left( 1 + \frac{\gamma_L}{s_E} \right), 
\end{equation}
and this holds trivially if $l \geq L-1$ since $\gamma_L \leq 4-t$. For $\varepsilon \ll s_E-1$ \eqref{lbound} gives that $l/L$ is bounded away from zero for $t>2$, with implicit constant depending only on $t$. Quantitatively:
\begin{equation} \label{fudgefactor2} \frac{s_E-1}{2} \leq \left(\frac{s_E-1}{2}\right)\left( 1 + \frac{\gamma_L}{s_E} \right) -O(\varepsilon) \leq \frac{l}{L} < 1. \end{equation} It follows that if $l \neq L$ (note $l \neq L-1$ always holds),
\begin{equation} \label{fudgefactor}  \gamma_l -\gamma_{l+1} \lesssim \frac{L(s_E-1)}{l(l+1)} \lesssim \frac{1}{\ell+1} \lesssim \frac{1}{L(s_E-1)} \lesssim \frac{\varepsilon}{s_E-1}. \end{equation}
 The left-hand side of \eqref{reusedbelow12} is trivially bounded by 
\begin{equation} \label{reuse} \int_0^{\pi} \int 
 \left(\delta^{-1} \left(P_{\mathbb{V}_{\theta \sharp}^{\perp}}\mu \right)\left\{ y \in B_E\left(x,\Delta_l \right) : d_E(\pi(x), \pi(y) )< \delta \right\}\right)^{p-1}  \\
 d\left( P_{\mathbb{V}_{\theta \sharp}^{\perp}}\mu \right)(x) \, d\theta. \end{equation}
The argument to bound \eqref{reuse} is almost identical to the argument bounding the left-hand side of \eqref{reusedbelow11} up until \eqref{reusedagain}, with $ 2^{-k}$ in \eqref{reusedbelow11} replaced by $\Delta_l$ in \eqref{reuse}, $\mu_{F_k}$ replaced by $\mu$, the condition $B \in \mathcal{B}_b$ defined to match the inequality in \eqref{zset2} rather than \eqref{zset}, and the set $\mathcal{B}$ is taken as $\mathcal{D}_{\Delta_l}(\supp \mu)$. If the terms from $\mathcal{B}_b$ dominate, or if the terms from $\mathcal{B}_g$ and $\Theta_b$ dominate, then the proof is identical as in those parts of the argument the lower density assumption is not used. Analogously to \eqref{referencedagain} and \eqref{intheabove}, the bound obtained in these cases will be much better than the sharp bound $2^{-k(p-1)(t-2)}$ provided that $p$ is taken sufficiently close to 1. 

If the low frequency case from $\mathcal{B}_g$ and $\Theta_g$ dominates, then by uniformity the bound analogous to \eqref{referencedagain2} in this case matches the bound of \eqref{pause33} below if only the $j = \log_2(\Delta_L)$ summand in \eqref{pause33} is kept. Therefore, the bound obtained in the low frequency case is at least as good as the bound obtained in the high frequency case, so it remains to check this final case. 

In the final case where terms from $\mathcal{B}_g$ and $\Theta_g$ dominate, and in the the high frequency case at the point analogous to \eqref{reusedagain}, applying \eqref{thirdestimate} in place of \eqref{secondestimate} to the second factor yields, in place of \eqref{reusedbelow21},
\begin{multline} \label{pause33} \eqref{reusedbelow12} \lesssim \sum_{B \in \mathcal{D}_{\Delta_{l}}(\supp \mu)_g} \left( \int_{\Theta_{g,B}} \mu(T_{B, \theta}) \, d\theta \right)^{1/q} \times \\
 \mu(B)^{1/q'} \left(\sum_{j =\left\lvert \log_2(\Delta_l)\right\rvert}^{\infty}  2^{j(1+\varepsilon)/q}  \sup_{D \in \mathcal{D}_{2^{-j}}(\supp \mu) } \mu(D)^{1/q}\right), \end{multline}
where each $T_{B, \theta}$ satisfies, analogously to \eqref{bushbound1} but using the more general term from \eqref{zset2}, 
\[ \mu(T_{B,\theta}) \lesssim \Delta_l^{2-\varepsilon}\sup_{r>0} \min\left\{ r \Delta_l^{-3} \sup_{D \in \mathcal{D}_{\Delta_l}(\supp \mu)} \mu(D), r^{t-3} \right\}. \]
This simplifies to 
\begin{equation} \label{goodtubebound} \mu(T_{B,\theta}) \lesssim 
\Delta_l^{2-\varepsilon} \left( \frac{ \left\lvert \mathcal{D}_{\Delta_l}(\supp \mu ) \right\rvert}{\Delta_l^{-3}}\right)^{-1+\frac{1}{4-t}}. \end{equation}
The second factor of each term in \eqref{pause33} satisfies 
\begin{multline} \label{keepstored} \mu(B)^{1/q'}\sum_{j =\left\lvert \log_2(\Delta_l)\right\rvert}^{\infty}  2^{j(1+\varepsilon)/q} \sup_{D \in \mathcal{D}_{2^{-j}}(\supp \mu) } \mu(D)^{1/q} \\
\lesssim \mu(B)^{1/q'} \sum_{m =l}^{L-1} \sum_{j = \left\lvert \log_2(\Delta_m)\right\rvert}^{\left\lvert \log_2(\Delta_{m+1})\right\rvert}  \Delta_{m+1}^{-(1+\varepsilon)/q} \sup_{D \in \mathcal{D}_{\Delta_m }(\supp \mu) } \mu(D)^{1/q} \\ +  \mu(B)^{1/q'} \sum_{j = \left\lvert \log_2(\Delta_{L})\right\rvert }^{\infty} 2^{j(1+\varepsilon)/q}  \min\left\{ 2^{-js_E/q}, \left\lvert \mathcal{D}_{\Delta_L}(\supp \mu) \right\rvert^{-1/q} \right\},    \end{multline}
By the maximality of $l$, and the uniformity of $\mu$ pigeonholed at the start of the proof, 
\begin{multline*} \eqref{keepstored} \lesssim \mu(B)^{1/q'}\Delta_l^{-O(\varepsilon/q)} \sum_{m =l}^{L-1} \sum_{j =\left\lvert \log_2(\Delta_m)\right\rvert}^{\left\lvert \log_2(\Delta_{m+1})\right\rvert}  \left( \Delta_{m+1}^{-1} \Delta_m^{s_E+\gamma_m} \right)^{1/q} \\ +  \mu(B)^{1/q'} \Delta_l^{-O(\varepsilon/q)} \left\lvert \mathcal{D}_{\Delta_L}(\supp \mu) \right\rvert^{-\left( 1- \frac{1}{s_E} \right)/q},  \end{multline*}
where the factor $\Delta_l^{-O(\varepsilon/q)}$ is needed to deal with $m =l$ term; the bound on $\mu(D)$ for $D \in  \mathcal{D}_{\Delta_l }(\supp \mu)$ is obtained by covering $D$ with $\Delta^{-O(1)}$ many $\Delta_{l+1}$ cubes and using \eqref{fudgefactor} and \eqref{fudgefactor2}.  By simplifying, this becomes
\begin{multline} \label{abovebecomes} \eqref{keepstored} \lesssim \mu(B)^{1/q'}\Delta_l^{-O(\varepsilon/q)} \sum_{m =l}^{L-1}  \left( \Delta_{m+1}^{-1} \Delta_m^{s_E+\gamma_m} \right)^{1/q}  \\
+  \mu(B)^{1/q'}\Delta_l^{-O(\varepsilon/q)}  \left\lvert \mathcal{D}_{\Delta_L}(\supp \mu) \right\rvert^{-\left( 1- \frac{1}{s_E} \right)/q}. \end{multline}
If $l < L$, then by the definition \eqref{gammaL1defn} of $\gamma_{L-1}$ the last term is bounded by the term with $m= L-1$: 
\[ \left\lvert \mathcal{D}_{\Delta_L}(\supp \mu) \right\rvert^{-\left( 1- \frac{1}{s_E} \right)/q} \leq    \Delta_L^{(s_E + \gamma_L)\left( 1- \frac{1}{s_E} \right)/q} 
\leq \Delta_l^{-O(\varepsilon/q)} \left( \Delta_{L}^{-1} \Delta_{L-1}^{s_E + \gamma_{L-1}}\right)^{1/q},  \]
so \eqref{abovebecomes} becomes
\begin{equation} \label{pause44}  \eqref{keepstored} \lesssim \mu(B)^{1/q'}   \Delta_l^{-O(\varepsilon/q)} \begin{cases} \sum_{m =l}^{L-1}  \left( \Delta_{m+1}^{-1} \Delta_m^{s_E+\gamma_m} \right)^{1/q} & l < L \\
\left\lvert \mathcal{D}_{\Delta_L}(\supp \mu) \right\rvert^{-\left( 1- \frac{1}{s_E} \right)/q} & l =L. \end{cases} \end{equation}
It will be shown that the terms in the sum are equal to each other:
\begin{equation} \label{constancy} \Delta_{m+1}^{-1} \Delta_m^{s_E+\gamma_m} = \Delta_{m+2}^{-1} \Delta_{m+1}^{s_E+\gamma_{m+1}}, \qquad l \leq m \leq L-2. \end{equation}
Since $\Delta_m = \Delta^{m}$, the above is equivalent to 
\[ -(m+1)+m(s_E + \gamma_m) = -(m+2) + (m+1)(s_E + \gamma_{m+1}), \qquad l \leq m \leq L-2, \]
which simplifies to 
\[ m \gamma_m = (m+1) \gamma_{m+1} + s_E-1, \qquad l \leq m \leq L-2, \]
which is equivalent to the defining recursive relation in \eqref{recursion}, so this establishes \eqref{constancy}. Using \eqref{constancy} to replace every term in \eqref{pause44} by the term with $m=l$ gives 
\[ \eqref{keepstored} \lesssim \mu(B)^{1/q'}\Delta_l^{-O(\varepsilon/q)}  \begin{cases}   \left( \Delta_{l+1}^{-1} \Delta_l^{s_E+\gamma_l} \right)^{1/q} & l < L \\
\left\lvert \mathcal{D}_{\Delta_L}(\supp \mu) \right\rvert^{-\left( 1- \frac{1}{s_E} \right)/q} & l =L. \end{cases}. \]
Using this to bound the second factor in \eqref{pause33} gives
\begin{multline} \label{pause34} \eqref{reusedbelow12} \lesssim \Delta_l^{-O(\varepsilon/q)}  \sum_{B \in \mathcal{D}_{\Delta_{l}}(\supp \mu)_g} \left( \int_{\Theta_{g,B}} \mu(T_{B, \theta}) \, d\theta \right)^{1/q}   \mu(B)^{1/q'} \times \\
\begin{cases}    \Delta_l^{(s_E+\gamma_l-1)/q} & l < L \\
\left\lvert \mathcal{D}_{\Delta_L}(\supp \mu) \right\rvert^{-\left( 1- \frac{1}{s_E} \right)/q} & l =L. \end{cases} \end{multline}
Applying \eqref{goodtubebound} to the first factor in \eqref{pause34}, and then using Hölder's inequality gives
\begin{multline*} \eqref{reusedbelow12} \lesssim \Delta_l^{-O(\varepsilon/q)} \left\lvert \mathcal{D}_{\Delta_{l}}(\supp \mu)\right\rvert^{1/q}  \left( \Delta_l^{2-\varepsilon} \left( \frac{ \left\lvert \mathcal{D}_{\Delta_l}(\supp \mu ) \right\rvert}{\Delta_l^{-3}}\right)^{-1+\frac{1}{4-t}} \right)^{1/q} \\
\times \begin{cases}    \Delta_l^{(s_E+\gamma_l-1)/q} & l < L \\
\left\lvert \mathcal{D}_{\Delta_L}(\supp \mu) \right\rvert^{-\left( 1- \frac{1}{s_E} \right)/q} & l =L. \end{cases} \end{multline*}
But $\left\lvert \mathcal{D}_{\Delta_{l}}(\supp \mu)\right\rvert \leq \Delta_l^{-(s_E+\gamma_l)}$ by the definition of $l$, and $\frac{1}{t-1}+\frac{1}{4-t} -1 >0$ when $2 < t < 3$, so this becomes
\[ \eqref{reusedbelow12} \lesssim \Delta_l^{-O(\varepsilon/q)}   \begin{cases}\Delta_l^{\left(t-2 + \gamma_l- \frac{\gamma_l}{4-t}\right)/q} & l < L \\
\Delta_L^{\left(t-2-\gamma_L\left(  \frac{1}{4-t}+ \frac{1}{t-1}-1\right)\right)/q} & l=L. \end{cases} \]
Since $\Delta_l = \Delta^l$, this yields 
\begin{equation} \label{pause67} \eqref{reusedbelow12} \lesssim 2^{kO(\varepsilon/q)}\begin{cases}  \Delta^{\left(l(t-2) + \frac{l \gamma_l(3-t)}{4-t}\right)/q} & l < L \\
\Delta_L^{\left(t-2-\gamma_L\left(  \frac{1}{4-t}+ \frac{1}{t-1}-1\right)\right)/q} & l=L. \end{cases} \end{equation}
But by the definition \eqref{gammadefn} of $\gamma_l$, and some algebra,
\[ l(t-2) + \frac{l \gamma_l(3-t)}{4-t} = (t-2)\left[ \frac{l}{4-t} + L \left( \frac{3-t}{4-t}\right)\left( 1+ \frac{\gamma_L}{t-1}\right)\right]. \]
If \eqref{fudgefactor2} is applied to the right-hand side of the above, the term $4-t$ cancels, resulting in 
\[ l(t-2) + \frac{l \gamma_l(3-t)}{4-t} \geq L(t-2)\left( 1+ \frac{\gamma_L}{t-1}\right) \left[ \frac{1}{2} -O(\varepsilon) \right]. \]
Substituting into \eqref{pause67} yields
\begin{equation}  \label{cases} \eqref{reusedbelow12} \lesssim 2^{kO(\varepsilon/q)} \begin{cases} 2^{-k\left(\frac{t-2}{2}\right)\left( 1+\frac{\gamma_L}{t-1} \right)/q} & l < L \\
2^{-k\left(t-2-\gamma_L\left(  \frac{1}{4-t}+ \frac{1}{t-1}-1\right)\right)/q} & l=L. \end{cases} \end{equation}

When $t>5/2$, a slightly better bound will be obtained when $l=L$ (still in the case where terms from $\mathcal{B}_g$ and $\Theta_g$ dominate). The high-low decomposition used is slightly different:
\[  \mu  = \mu \ast \widecheck{\psi_{l}} + \mu \ast \widecheck{\psi_h}, \]
where $\psi_l$ is a smooth bump on $|\xi| \lesssim 2^j$, and $\psi_h$ is a smooth bump on $|\xi| \gtrsim 2^j$,  where
\begin{equation} \label{jformula} j := k\left( 1 + \frac{\gamma_L}{4-t}\right). \end{equation}
Since $t>5/2$, this $j$ satisfies $j \geq \frac{s_E+\gamma_L}{s_E}k$. Define $q$ by $\frac{1}{q} = p-1$. Repeating a similar argument to the above in the high case, the analogous bound to \eqref{pause34} is the simpler bound
\[ \eqref{reusedbelow12} \lesssim \Delta_L^{-O(\varepsilon/q)}  \sum_{B \in \mathcal{D}_{\Delta_{L}}(\supp \mu)_g} \left( \int_{\Theta_{g,B}} \mu(T_{B, \theta}) \, d\theta \right)^{1/q}   \mu(B)^{1/q'} 2^{-j(s_E-1-\epsilon)/q}. \]
Using \eqref{goodtubebound}, Hölder's inequality, and $\left\lvert \mathcal{D}_{\Delta_L}( \supp \mu) \right\rvert \leq \Delta_L^{-(s_E+\gamma_L) }$, this simplifies to
\begin{equation} \label{highcasebound} \eqref{reusedbelow12} \lesssim 2^{\left[ k O(\varepsilon) + k\frac{\gamma_L}{4-t}  - j(t-2)\right]/q}. \end{equation}
In the low case, still with $l =L$, by Hölder's inequality the expression in \eqref{reusedbelow12} is bounded by 
\begin{multline} \label{youngandkakeya} 2^{-k/q} \\
\times \left(\int_0^{\pi} \int 
 \left(2^k\delta^{-1} \left(P_{\mathbb{V}_{\theta \sharp}^{\perp}}\mu \right)\left\{ y \in B_E\left(x,2^{-k}\right) : d_E(\pi(x), \pi(y) )< \delta \right\}\right) \, dx \, d\theta \right)^{1/q}  \\
\times \left( \int_0^{\pi} \int \left\lvert P_{\mathbb{V}_{\theta \sharp}^{\perp}}\left(\mu \ast \phi_j \right)(x)\right\rvert^{q'} \, dx \, d\theta\right)^{1/q'}, \end{multline} 
where $\frac{1}{q} = p-1$, and $\phi_j$ is a smooth bump function rapidly decaying outside on $B_E(0, 2^{-j})$, with $\|\phi_j\|_{\infty} \lesssim 2^{j(3+\varepsilon)}$. By applying Young's convolution inequality (or just Fubini) to the first factor in \eqref{youngandkakeya}, and applying \cite[Theorem~3.1]{harris2023} (the $L^{q'}$ version that interpolates between the $L^{3/2}$ bound and the trivial $L^1$ bound) to the second factor, this gives
\begin{equation} \label{pause337} \eqref{reusedbelow12} \lesssim 2^{-k/q} c_{3+\varepsilon^2}( \mu \ast \phi_j )^{1/q}.\end{equation}
By Lemma~\ref{euclidean} and \eqref{comparisonfrostman}, 
\[ c_{3+\varepsilon^2}(\mu \ast \phi_j) \lesssim 2^{kO(\varepsilon)} \sup_{r>0} \min\left\{ r2^{j(4-t)}, r^{t-3} \right\} = 2^{kO(\varepsilon)} 2^{j(3-t)}, \]
so \eqref{pause337} becomes
\[ \eqref{reusedbelow12} \lesssim 2^{\left( -k(1+O(\varepsilon)) +j(3-t)\right)/q}. \]
The formula \eqref{jformula} implies that the above bound in the low case matches the bound \eqref{highcasebound} in the high case, and hence
\[ \eqref{reusedbelow12} \lesssim  2^{-k\left[  t-2  - \frac{\gamma(3-t)}{4-t} \right]/q}. \]
Combining with \eqref{cases} yields 
\begin{equation} \label{cases2} \eqref{reusedbelow12} \lesssim 2^{kO(\varepsilon/q)} \begin{cases} 2^{-k\left(\frac{t-2}{2}\right)\left( 1+\frac{\gamma_L}{t-1} \right)/q} & l < L \\
2^{-k\left(t-2-\gamma_L\min\left\{ \frac{1}{4-t}+ \frac{1}{t-1}-1, \frac{3-t}{4-t} \right\}\right)/q} & l=L. \end{cases} \end{equation}
But $1/q = p-1$ and by the definition \eqref{gammaLdefn2} of $\gamma_L$, the two bounds are the same when $\gamma_L < 4-t$. If $\gamma_L = 4-t$ then $l=L$ since $s_E + \gamma_L = 3$, so in any case the equation for $\gamma_L$ can be substituted into the second formula in \eqref{cases2}. When $5/2 \leq t < 3$ this gives
\[ \eqref{reusedbelow12} \lesssim 2^{kO(\varepsilon/q)}2^{-k\left(t-2-\gamma_L\left(\frac{3-t}{4-t}\right)\right)/q}.  \]
Substituting the value of $\gamma_L$ from \eqref{gammaLdefn} in the case $5/2 \leq t < 3$ verifies \eqref{reusedbelow12} when $5/2 \leq t < 3$. If $2 < t < 5/2$, then $\gamma_L < 4-t$, so the first formula in \eqref{cases2} holds as both bounds are the same, which yields
\[ \eqref{reusedbelow12} \lesssim 2^{kO(\varepsilon/q)}2^{-k\left(\frac{t-2}{2}\right)\left( 1+\frac{\gamma_L}{t-1} \right)/q}.  \]
Substituting the value of $\gamma_L$ from \eqref{gammaLdefn} in the case $2 < t < 5/2$ verifies \eqref{reusedbelow12} when $2 < t < 5/2$. \end{proof}
\begin{remark} \label{gammaLremark} Taking $\gamma_L=0$ in the above proof instead of \eqref{gammaLdefn} would only yield the weaker lower bound of $1+\frac{1}{2} \dim A$ near 2, in place of the bound in Theorem~\ref{etccorollary2}.  \end{remark}

\bibliographystyle{alpha}
\bibliography{main}

\end{document}